\documentclass[10pt]{article}

\usepackage[usenames,dvipsnames]{pstricks}
\usepackage{epsfig}
\usepackage{pst-grad} % For gradients
\usepackage[space]{grffile} % For spaces in paths
\usepackage{etoolbox} % For spaces in paths
\makeatletter % For spaces in paths
\patchcmd\Gread@eps{\@inputcheck#1 }{\@inputcheck"#1"\relax}{}{}
\makeatother

\usepackage{graphicx}
\usepackage{amsmath}
\usepackage{psfrag}
\usepackage{color,soul}
\usepackage{natbib}
\usepackage{siunitx}
\usepackage{float}
\usepackage{tikz}
\usepackage{multirow}
\usepackage{bm}
\usepackage{textgreek}
\usepackage{url}
\usepackage{amssymb}
\usepackage{wasysym}
\usepackage{xfrac}
\usepackage{setspace}
\usepackage{gensymb}
\usepackage{mathtools}
\usepackage{arydshln}
\usepackage[percent]{overpic}
\usepackage{xspace}
\usepackage{bold-extra}
\usepackage[most]{tcolorbox}
\usepackage{subcaption}
\usepackage{array}
\usepackage{makecell}
\usepackage{xcolor}
\usepackage{pict2e}
\usepackage{bigints}
\usepackage{array}
\usepackage{amsopn}
\newcolumntype{P}[1]{>{\centering\arraybackslash}p{#1}}

\DeclareMathOperator*{\argmax}{arg\,max}
\definecolor{myred}{rgb}{0.5, 0, 0.0}
\definecolor{myblue}{rgb}{0, 0, 0.5}

\usepackage[noblocks]{authblk}
\usepackage[margin=1in]{geometry}
\usepackage[title]{appendix}

\newcommand\affiliation[1]{\gdef\@affiliation{\let\aff\aff@inst#1}}
\gdef\@affiliation{}
\def\email#1{Email address for correspondence: #1}
\def\aff#1{\ignorespaces\textsuperscript{#1}}
\def\corresp#1{\unskip\thanks{#1}}
\numberwithin{equation}{section}

\renewenvironment{abstract}
{\begin{quote}
\noindent \rule{\linewidth}{.5pt}\par{\bfseries \abstractname.}}
{\medskip\noindent \rule{\linewidth}{.5pt}
\end{quote}
}

  \DeclareTextFontCommand\textsfi{\usefont{OT1}{cmss}{m}{sl}}
  \DeclareMathAlphabet\mathsfi            {OT1}{cmss}{m}{sl}
  \DeclareTextFontCommand\textsfb{\usefont{OT1}{cmss}{bx}{n}}
  \DeclareMathAlphabet\mathsfb            {OT1}{cmss}{bx}{n}
  \DeclareTextFontCommand\textsfbi{\usefont{OT1}{cmss}{m}{sl}}
  \DeclareMathAlphabet\mathsfbi            {OT1}{cmss}{m}{sl}
  \IfFileExists{t1phv.fd}
  {\DeclareTextFontCommand\textsfbi{\usefont{T1}{phv}{b}{it}}
  \DeclareMathAlphabet\mathsfbi            {T1}{phv}{b}{it}
  }{}
  \IfFileExists{ot1phv.fd}
  {\DeclareTextFontCommand\textsfbi{\usefont{OT1}{phv}{b}{it}}
  \DeclareMathAlphabet\mathsfbi            {OT1}{phv}{b}{it}
  }{}
  
\DeclareSymbolFont{matha}{OML}{txmi}{m}{it}% txfonts

\usepackage{hyperref}

\definecolor{darkblue}{rgb}{0,0,0.80}

\hypersetup{
    colorlinks=true,
    linkcolor={darkblue},
    citecolor={darkblue},
    filecolor=black,
    urlcolor={darkblue},
    plainpages=false,
}

\title{\bf Space-time POD and the Hankel matrix}

\author[1]{\bf Peter Frame\corresp{\email{pframe@umich.edu}}}
\author[1]{\bf Aaron Towne}

\affil[1]{\normalsize Department of Mechanical Engineering, University of Michigan, Ann Arbor, MI, USA \vspace{-1cm}}

\date{}
\begin{document}
\maketitle

% Abstract
\begin{abstract}
Time-delay embedding is an increasingly popular starting point for data-driven reduced-order modeling efforts. In particular, the singular value decomposition (SVD) of a block Hankel matrix formed from successive delay embeddings of the state of a dynamical system lies at the heart of several popular reduced-order modeling methods. In this paper, we show that the left singular vectors of this Hankel matrix are a discrete approximation of classical space-time proper orthogonal decomposition (POD) modes, and the singular values are square roots of the POD energies. This connection establishes a clear interpretation of the Hankel modes grounded in classical theory, and we gain insights into the Hankel modes by instead analyzing the equivalent discrete space-time POD modes in terms of the correlation matrix formed by multiplying the Hankel matrix by its conjugate transpose. These insights include the distinct meaning of rows and columns, the implied norm in which the modes are optimal, the impact of the time step between snapshots on the modes, and an interpretation of the embedding dimension/height of the Hankel matrix in terms of the time window on which the modes are optimal. Moreover, the connections we establish offer opportunities to improve the convergence and computation time in certain practical cases, and to improve the accuracy of the modes with the same data. Finally, popular variants of POD, namely the standard space-only POD and spectral POD, are recovered in the limits that snapshots used to form each column of the Hankel matrix represent flow evolution over short and long times, respectively. \\  % DO NOT DELETE THE \\
\end{abstract}

%%%%%%%%%%%%%%%%%%%%%%%%%%%%%%%%%%%%%%%%%%%%%%%%%%%%%%%%%%%%%%%%%%%%%%%%%%%%%%%
% -----------------------------------------------------------------------------
% -- Introduction
% -----------------------------------------------------------------------------

\section{Introduction}
\label{sec:intro}

Time series data, generated from simulations or experiments, are abundant in science and engineering, but analyzing or interpreting these data can be challenging. Often, such as in climate science or the analysis of financial data, researchers want an understanding of the governing laws underlying the time series. Other times, as in the case of fluid mechanics, a precise physical model exists, but it may be difficult to interpret the results, or simulating the full model may be computationally costly. In such cases, the goal is often to use data to uncover key physical mechanisms that contribute to the underlying dynamics \cite{Rowley00} or to derive a less physical, but more computationally efficient, reduced-order model capable of approximating the dynamics at low cost.
\\

Many techniques have emerged to address these challenges. At the heart of several of them, especially ones popular in the dynamical systems community, is analysis of a Hankel matrix. A Hankel matrix has constant skew diagonals, i.e., the ($i,j$) entry of the matrix only depends on $i+j$. Common applications include system identification and minimal system realization \cite{Fazel}, and the use of the Hankel matrix in this context goes back to the 1960s \cite{Silverman}. The Hankel matrix plays a key role in singular spectral analysis (SSA) \cite{VAUTARD89} and in the eigensystem realization algorithm (ERA) \cite{Juang85}. Recently, the Hankel matrix has been used in the context of fluid dynamics. It is central to balanced truncation \cite{Moore81}, which was made scalable to fluid dynamics problems with balanced POD \cite{Willcox02,Rowley05}. The Hankel matrix is also used in a variant of dynamic mode decomposition (DMD) \cite{Schmid10} called Hankel DMD \cite{Arbabi}, where the goal is to extract the spectrum and modes of the Koopman operator for some dynamical system by performing DMD on a Hankel matrix of observables. 
\\

The entries of the Hankel matrix are taken from time series data, and moving down or right along the columns or rows of the matrix corresponds to moving forward in time. Therefore, the columns of the Hankel matrix are delay embeddings of the dynamical system that produced the time series. Delay embedding, which goes back to work from the 1980s \cite{Packard80,Takens81}, is a method of encoding the state of a dynamical system by recording the time history of one (or a few) of its observables. Intuitively, the state of an $n$-dimensional system should be determined by $n$ independent observables, which don’t need to be the original degrees of freedom of the system \cite{Packard80}. This idea was later made rigorous by Takens \cite{Takens81}, and, under weak conditions, $2n+1$ entries of a time series are needed to determine the state of the system.
\\

A recurring theme in applications of the Hankel matrix in dynamical systems is the singular value decomposition (SVD). Both SSA and ERA obtain their bases from the SVD of the Hankel matrix. In balanced truncation and balanced POD, the Hankel singular values are used. In the Hankel alternative view of Koopman (HAVOK) framework \cite{Brunton16}, the SVD of a Hankel matrix of data is used in order to form a low-rank linear model of the dynamics on some chaotic attractor. The modes in their linear model are the left singular vectors of the Hankel matrix formed from a time series of the dynamical system. Convolutional coordinates, used to represent the state of a dynamical system at a particular time in terms of its representation in some predefined temporal basis, can be defined in terms of the continuous SVD (Schmidt decomposition) of a continuous Hankel matrix \cite{Kamb20}. The left singular vectors of the Hankel matrix have recently been called principal component trajectories and used for control \cite{Dylewsky22}.
\\

A second, older, technique for analyzing time series data is proper orthogonal decomposition (POD). Originally introduced to the fluid dynamics community by Lumley \cite{Lumley67} in 1967, it is known by a variety of names in other areas including principal component analysis, Karhunen-Lo\`eve decomposition, and empirical orthogonal functions. In POD, the flow data is analyzed statistically, and the objective is to search for the modes that most efficiently represent the data. Specifically, POD modes are defined to minimize the reconstruction error, as measured by the average square inner product, compared to any other basis of the same dimension. As introduced by Lumley, the most general version of POD seeks to describe the time evolution of the flow for a prespecified window of time, so the basis functions are functions both of space and of time. This general version is called space-time POD, and the reconstruction of a flow over the time window consists of these basis functions multiplied by constant coefficients. Space-time POD has been used infrequently in the literature. Notable exceptions include application of space-time POD to optimally describe transients \cite{Gordeyev13}, generalize dimension reduction methods \cite{delRosario18}, and study acoustic intermittency in the form of bursts in jets via conditioning \cite{Schmidt18CTR,Schmidt19}.
\\

Following the work of Sirovich \cite{Sirovich87} and Aubry \cite{Aubry91}, today the most popular form of POD involves modes that are functions of space only. To represent a time-dependent flow, these modes are multiplied by time-varying coefficients. We refer to this form of POD as space-only POD \cite{Aaron18}, but note that it is often referred to simply as POD in the literature. Space-only POD has been used extensively to form Galerkin-based reduced-order models \cite{Aubry88,Rowley04,Rowley17}, educe physically meaningful structures from flow data \cite{Holmes12,Moin89}, and reduce the data needed to store flow data \cite{Pollard16}. We will show that space-only POD can be understood within the more general space-time POD framework as the limit as the time interval on which the optimization problem is defined goes to zero.
\\

An increasingly popular variant of POD is spectral POD (SPOD). While also introduced in the original work of Lumley \cite{Lumley67,Lumley70}, this form of POD was rarely employed until recently \cite{Aaron18,Cavalieri19,Symon21,Schmidt18}. Here, the objective is to optimally represent statistically stationary flows in the frequency domain. At each frequency, SPOD provides a set of spatial modes that capture the portion of the flow at that frequency more accurately, on average, than any other basis of the same order. SPOD can be formulated as the limit of space-time POD as the time interval on which the modes are defined goes to infinity. 
\\

Our aim in this paper is to show that the singular modes of the Hankel matrix, i.e., principal component trajectories, are a discrete approximation of the classical space-time POD modes and to demonstrate that understanding them as such is useful for analyzing and improving their properties. This connection is established by observing that the Hankel matrix multiplied by its Hermetian transpose provides an approximation of the space-time correlation matrix whose eigendecomposition defines discrete space-time POD modes. We show that analyzing the space-time correlation matrix that would be formed using a particular Hankel matrix leads to insight into the modes and, in some cases, guidance about how to improve them. The Hankel modes are not the only approximation one could form of the space-time POD modes using the available data, and we show that, in several cases, they are not the most practical approximation. With a surplus of underlying data, we show that a sufficiently accurate correlation can be formed by throwing out many of the columns in the Hankel matrix, which can drastically reduce the computational cost of the SVD. We also show that a more accurate approximation of the correlation can be obtained by fully exploiting the ergodicity of the system (ergodicity is assumed in forming the Hankel matrix). For low dimensional systems, we show that if little data is available, this improves accuracy significantly, and that if much data is available, this method is asymptotically faster than the Hankel matrix approach. The connection to space-time POD also makes clear the assumed inner product, which defines the sense in which the Hankel modes are optimal, and clarifies the impact of the time step between successive snapshots on the approximation. Additionally, we show that the height and width of the Hankel matrix determine the extent to which temporal correlation is accounted for in defining the modes, and the convergence of the correlation, respectively. Finally, we show that in the limits of short and long delays, space-only and spectral POD are recovered, respectively.
\\

The remainder of the paper is organized as follows. In section \ref{sec:Hankel Intro}, we define the Hankel matrix and reiterate some of its applications. In section \ref{sec:SOPOD}, we motivate and derive the continuous and discrete forms of space-only POD. We do this using a formalism that makes spectral and space-time POD follow easily, and that highlights the central role of the correlation matrix. In section \ref{sec:SPOD}, we derive the continuous and discrete forms of spectral POD. In section \ref{sec:STPOD}, we derive the continuous form of space-time POD and show that, in the limits of short and long time intervals, it reduces to space-only and spectral POD, respectively. In section \ref{sec:Hankel}, we show that the singular modes of the Hankel matrix provide a discrete approximation of space-time POD modes by noting the Hankel matrix can be used to approximate the space-time correlation tensor. This connection motivates several improvements to the Hankel SVD procedure, and we also prove results about the convergence to space-only and spectral POD analogous to those for the continuous case. In section \ref{sec:exploiting_ergodicity}, we show that with the same time series data, we can construct a more accurate correlation matrix than the Hankel matrix multiplied by its conjugate transpose by exploiting the ergodicity of the system. The eigendecomosition of this correlation matrix gives more converged modes than the SVD of the Hankel matrix. In section \ref{sec:example}, we demonstrate these results using a lid-driven cavity flow at $Re = 22,000$. Finally, in section \ref{sec:conclusions}, we summarize the paper and present our conclusions.

\section{The Hankel matrix} \label{sec:Hankel Intro}
A Hankel matrix is a matrix whose skew-diagonals are constant, i.e., the ($i,j$) entry only depends on $i+j$,
\begin{equation} \label{HankelDef}
    H = \begin{bmatrix} 
    g_1 & g_2 & \dots & g_m\\ g_2 & g_3  & \dots & g_{m+1}\\ 
    \vdots & \vdots & \ddots & \vdots\\
    g_d &  g_{d+1}  & \dots    & g_{m+d-1} 
    \end{bmatrix}
    \text{.}
\end{equation}
In dynamical systems theory, these matrices are formed from a time series of observables $g_1, g_2, g_3, \dots ,g_{m+d-1}$. This means that the columns of $H$ are time delay embeddings of the dynamical system. In this paper, we focus our attention on the more general block Hankel matrix
\begin{equation}
   {\bf{H}} =  \begin{bmatrix} {\bf{q}}_1 & {\bf{q}}_2 & \dots & {\bf{q}}_m \\ {\bf{q}}_2 & {\bf{q}}_3  & \dots & {\bf{q}}_{m+1}\\ 
    \vdots & \vdots & \ddots & \vdots\\
    {\bf{q}}_d &  {\bf{q}}_{d+1}  & \dots  & {\bf{q}}_{m+d-1} 
    \end{bmatrix}\text{,}
\end{equation}
formed from a time series of vector-valued observables ${\bf{q}}_1,  {\bf{q}}_2, \dots , {\bf{q}}_d, \dots, {\bf{q}}_{m+d-1}$.
\\

Applications of the Hankel matrix include the singular spectrum analysis (SSA) \cite{VAUTARD89}, the eigensystem realization algorithm (ERA) \cite{Juang85} and the Hankel alternative view of Koopman (HAVOK) \cite{Brunton16}. All of these methods depend on the SVD of the Hankel matrix,
\begin{equation}
    {\bf{H}} = {\bf{U}}{\bf{\Sigma}}{\bf{V}}^* \text{,}
\end{equation}
and in particular use the left singular vectors of the Hankel matrix as a basis, which have been called principal component trajectories \cite{Dylewsky22}. The connections to time-delay embedding, Koopman theory, and dynamic mode decomposition have garnered increasing interest in analysis of the Hankel matrix.

\section{Space-only POD} \label{sec:SOPOD}
All forms of POD are statistical methods and view the dynamical system to which they are applied as random. This is a practical choice; though the system is not random, it might be chaotic, and without knowledge of the exact initial condition of the system, viewing it as random is the best we can do. If we think of a spatial realization of the flow, or a \textit{snapshot} $\boldsymbol{q}({\bf{x}})$, as a vector in a vector space, space-only POD seeks to find the direction in the vector space along which there is the most variation between different snapshots, i.e., the most energy as defined by some spatial inner product. Knowing the coordinate of a snapshot of the flow along this direction provides a better approximation, on average, than any other coordinate. To characterize a flow exactly, we need to specify the flow field everywhere in the domain, but if giving the flow's coordinates along just a few important directions approximates it to high precision, then it is useful to look for these important directions. 
\\

\subsection{Formulation}
The direction with the most variation is formalized as the mode which maximizes the expected value of the square of the projection of the flow snapshot, $\boldsymbol{q}({\bf{x}})$ onto the mode, 
\begin{subequations} \label{SOPOD:Optimization}
\begin{equation} \label{Sop:maximization}
    \lambda[\boldsymbol{\phi}({\bf{x}})] =  \frac{1}{\| \boldsymbol{\phi}({\bf{x}})\| }  \mathbb{E}\big[ \|\langle \boldsymbol{q}({\bf{x}}), \boldsymbol{\phi} ({\bf{x}})\rangle \|^2 \big] \text{,}
\end{equation}
\begin{equation} \label{Sop:maximization_with_functional}
    \boldsymbol{\phi}_1({\bf{x}})  = \argmax \lambda[\boldsymbol{\phi}({\bf{x}})] \text{.}
\end{equation}
\end{subequations}
Here, the expectation operator $\mathbb{E}[\cdot ]$ acts over time, and $\lambda$ is a functional which takes any function as input (in the appropriate function space) and returns the variation in the flow along it. The first POD mode $\boldsymbol{\phi}_1$ maximizes this functional. Typically, the modes are normalized to unity, but for clarity we will keep the magnitude in the formulae. The magnitude operator $\| \cdot \|$ is defined in the usual way by the inner product, and the inner product takes the form 
\begin{equation} \label{eq:SOPOD_IP}
      \langle \boldsymbol{q}_1({\bf{x}}),\boldsymbol{q}_2({\bf{x}}) \rangle  = \int \limits_{\Omega}\boldsymbol{q}_2^{{\bf{*}}}({\bf{x}})\boldsymbol{W}({\bf{x}})\boldsymbol{q}_1({\bf{x}})d{\bf{x}} \text{,}
\end{equation}
where $\boldsymbol{W}$ is a weight matrix and $\Omega$ is the spatial domain of interest. The weight matrix is often used to make the norm correspond to some physical definition of energy, e.g., turbulent kinetic energy, or to give preference to certain flow variables or regions of the flow. The first POD mode $\boldsymbol{\phi}_1$ maximizes this functional.
In solving (\ref{Sop:maximization_with_functional}), it is helpful to rewrite (\ref{Sop:maximization}) as an inner product with the correlation tensor \cite{Lumley67,RowleyPhD},
\begin{equation} \label{PODCorrelation}
    \lambda[\boldsymbol{\phi}({\bf{x}})] =  \frac{1}{\|\boldsymbol{\phi}({\bf{x}})\|} \Big \langle \big \langle \boldsymbol{\phi}({\bf{x}}'),\boldsymbol{C}({\bf{x}},{\bf{x}}') \big \rangle,\boldsymbol{\phi}({\bf{x}}) \Big \rangle  \text{,}
\end{equation}
where the correlation is defined as
\begin{equation} \label{Sop: Correlation}
    \boldsymbol{C}({\bf{x}}_1,{\bf{x}}_2) = \mathbb{E}[\boldsymbol{q}({\bf{x}}_1) \boldsymbol{q}({\bf{x}}_2)^{{\bf{*}}}] \text{.}
\end{equation}
The correlation tensor is symmetric (under interchange of its two arguments), so it has orthogonal eigenfunctions $\{ \boldsymbol{\nu}_1({\bf{x}}), \boldsymbol{\nu}_2({\bf{x}}), \dots \}$, which can be ordered so that the associated eigenvalues $\{\lambda_1 \geq \lambda_2 \geq , \dots \geq 0\}$ are non-increasing. These are eigenfunctions in the sense that they satisfy
\begin{equation}
    \langle \boldsymbol{\nu}_k({\bf{x}}_2),\boldsymbol{C}({\bf{x}}_1,{\bf{x}}_2) \rangle = \lambda_k\boldsymbol{\nu}_k({\bf{x}}_1) \text{.}
\end{equation}
Writing $\boldsymbol{\phi}$ in the basis of these eigenfunctions provides insight into the maximization problem (\ref{SOPOD:Optimization}), and the energy of $\boldsymbol{\phi}$ (\ref{PODCorrelation}) can be written in terms of the expansion coefficients used to express it \cite{Lumley67},
\begin{equation} \label{sec:SOPOD:eq:Energy_corr_eig}
    \lambda[\boldsymbol{\phi}({\bf{x}})] = \sum_{k=1}^\infty |c_k|^2 \lambda_k  \quad \text{for normalized $\boldsymbol{\phi}({\bf{x}})$,}
\end{equation}
where the expansion coefficient $c_k = \langle \boldsymbol{\phi}(\boldsymbol{x}),\boldsymbol{\nu}_{k}(\boldsymbol{x}) \rangle$ is the projection of $\boldsymbol{\phi}({\bf{x}})$ on the $k^{\text{th}}$ eigenfunction. The solution that maximizes (\ref{sec:SOPOD:eq:Energy_corr_eig}) is $c_1 = 1, c_{\neq 1} = 0$, because $\lambda_1 > \lambda_{\neq 1}$. This tells us that the first POD mode is the eigenvector of the correlation tensor with the greatest eigenvalue. To define the latter modes, we maximize $\lambda[\boldsymbol{\phi}({\bf{x}})]$ over all $\boldsymbol{\phi}({\bf{x}})$ orthogonal to previous modes. A simple inductive argument shows that if the first $k$ POD modes are the first $k$ eigenfunctions of the correlation tensor, then the $k+1^{\text{st}}$ mode must be the $k+1^{\text{st}}$ eigenfunction because this orthogonality condition tells us that $c_{<k+1} = 0$ for the $k+1^{\text{st}}$ mode. The POD modes are therefore the eigenfunctions of the correlation tensor, $\boldsymbol{\phi}_k = \boldsymbol{\nu}_k$, and their energies are the eigenvalues, $\lambda[\boldsymbol{\phi}_k({\bf{x}})] = \lambda_k$. In terms of the inner product, the POD modes satisfy
\begin{equation} \label{SOPOD:integral_evp}
        \int \limits_{\Omega}\boldsymbol{C}({\bf{x}}_1,{\bf{x}}_2)\boldsymbol{W}({\bf{x}}_2)\boldsymbol{\phi}_{k}({\bf{x}}_2)d{\bf{x}}_2 = \lambda_k\boldsymbol{\phi}_{k}({\bf{x}}_1) \text{.}
\end{equation}
\subsection{With discrete data} \label{susbsec:SOPOD_disc}
In practice, space-only POD modes are approximated using data from a simulation or experiment defined on a discrete set of points and sampled in time. To approximate the true continuous space-only POD modes, the data is used to approximate a correlation matrix whose eigendecomposition gives the modes. Denoting one discrete snapshot as ${\bf{q}}_j \in \mathbb{R}^N$, we need a correlation matrix of the form ${\bf{C}}^{\bf{x}} \in \mathbb{R}^{N \times N}$, where each element represents the correlation between two components of the flow at different points. We label this spatial correlation matrix with an ${\bf{x}}$ superscript to distinguish it from other correlations that arise later in the paper. It can be approximated as
\begin{equation} \label{eqn:SOPOD_CorrelationTEns}
{\bf{C}}^{\bf{x}} = \frac{1}{m}{\bf{Q}}{\bf{Q}}^{{\bf{*}}} \quad \text{with} \quad {\bf{Q}}  = [{\bf{q}}_1,{\bf{q}}_2,\dots {\bf{q}}_m] \in \mathbb{R}^{N \times m} \text{.}
\end{equation}
The data matrix ${\bf{Q}}$ contains an ensemble of snapshots, and $\frac{1}{m}{\bf{Q}}{\bf{Q}}^{{\bf{*}}}$ is an approximation of ${\bf{C}}$ because it gives each component of the correlation as a sample average over the realizations,
\begin{equation}
    ({\bf{Q}}{\bf{Q}}^*)_{ij} = \sum_{k=1}^m ({\bf{q}}_{k})_i({\bf{q}}_{k}^*)_j \text{.}
\end{equation}
Discrete space-only POD modes are given by the eigenvectors of the discrete correlation matrix multiplied by the weight,
\begin{equation}\label{eq:SOPODEig1}
{\bf{C}}^{\bf{x}}{\bf{W}}{\bf{\Phi}} = {\bf{\Phi}}{\bf{\Lambda}}\text{,}
\end{equation}
where ${\bf{W}}$ is the discrete weight matrix. This weight is used both as the discrete version of the continuous weight and to account for numerical quadrature of the integral in (\ref{SOPOD:integral_evp}). The columns of the matrix ${\bf{\Phi}} = [{\bf{\phi}}_1, {\bf{\phi}}_2 \dots ]$ are the discrete space-only POD modes, and the diagonal matrix ${\bf{\Lambda}}$ contains the corresponding eigenvalues. 
\\

The discrete space-only POD modes are also related to the left singular vectors of the data matrix ${\bf{Q}}$ and the weight matrix, and the details of this relation will prove important to understanding the connection between the SVD of the Hankel matrix and POD. If we take the singular value decomposition to obtain 
\begin{equation} \label{eq:SOPOD:SVD}
\frac{1}{\sqrt{m}}{\bf{W}}^{\frac{1}{2}}{\bf{Q}} = {\bf{U}}{\bf{\Sigma}}{\bf{V}}^{{\bf{*}}} \text{,}
\end{equation}
then we can write the correlation tensor multiplied by the weight as ${\bf{C}}^{\bf{x}}{\bf{W}} = {\bf{W}}^{-\frac{1}{2}}{\bf{U}}{\bf{\Sigma}}{\bf{V}}^{{\bf{*}}}{\bf{V}}{\bf{\Sigma}}{\bf{U}}^{{\bf{*}}}{\bf{W}}^{\frac{1}{2}}$. Because both ${\bf{U}}$ and ${\bf{V}}$ are orthonormal, multiplying by ${\bf W}^{-\frac{1}{2}} {\bf U}$ we have 
\begin{equation}
{\bf{C}}^{\bf{x}}{\bf{W}}{\bf W}^{-\frac{1}{2}} {\bf U} = {\bf{W}}^{-\frac{1}{2}}{\bf{U}}{\bf{\Sigma}}^2 \text{.}
\end{equation}
This constitutes an eigendecomposition of ${\bf{C}}^{\bf{x}}{\bf{W}}$, and therefore the POD modes are related to the left singular vectors of the data matrix and the energies to the singular values,
\begin{equation}
{\bf{\Phi}} = {\bf{W}}^{-\frac{1}{2}}{\bf{U}} \quad \text{and} \quad {\bf{\Lambda}} = {\bf{\Sigma}}^2 \text{.}
\end{equation}
\\

It is important to remember that the correlation matrix is approximate because it is obtained from finite data. If the data are highly correlated, e.g., if the snapshots are taken from a time series whose length is on the same order as the characteristic timescale of the flow or shorter, then the correlation tensor, and hence the POD modes, will be inaccurate. Both will increase in accuracy with the number of realizations and the independence of the realizations.

\section{Spectral POD} \label{sec:SPOD}
Spectral POD produces an optimal frequency domain representation of statistically stationary flows. At each frequency, SPOD modes reconstruct the Fourier transform of a flow more accurately, on average, than any other reconstruction of the same order. 
\\

\subsection{Formulation}
Spectral POD can also be cast as an optimization problem, analogous to that of space-only POD, as 
\begin{subequations} \label{SPOD: optimization}
\begin{equation} \label{SPOD: maximization}
    \lambda_\omega[\boldsymbol{\phi}({\bf{x}})] =  \frac{1}{\| \boldsymbol{\phi}({\bf{x}})\| }  \mathbb{E}\big[ \|\langle \hat{\boldsymbol{q}}_{\omega}({\bf{x}}), \boldsymbol{\phi} ({\bf{x}})\rangle \|^2 \big] \text{,}
\end{equation}
\begin{equation}
    \boldsymbol{\psi}_{\omega,1}({\bf{x}})  = \argmax \lambda_{\omega}[\boldsymbol{\phi}({\bf{x}})] \text{.}
\end{equation}
\end{subequations}
Here, the expectation operator acts over the Fourier transform of segments of the flow. These segments can be subsections of a longer time series or correspond to separate realizations of the flow. Similar to the space-only case, $\lambda_{\omega}$ is a functional that returns the energy in the flow at frequency $\omega$ captured by the argument. The Fourier transformed flow is defined as 
\begin{equation}
    \hat{\boldsymbol{q}}_{\omega}({\bf{x}}) = \int \limits_{-\infty}^{\infty} \boldsymbol{q}({\bf{x}},t)e^{-i\omega t}dt \text{,}
\end{equation}
and the inner product is the same as before, given in (\ref{eq:SOPOD_IP}). Analogous to the space-only case, the optimization problem (\ref{SPOD: maximization}) can be rewritten in terms of correlations,
\begin{equation} \label{SPODCorrelation}
    \lambda_{\omega}[\boldsymbol{\phi}({\bf{x}})] =  \frac{1}{\|\boldsymbol{\phi}({\bf{x}})\|} \Big \langle \big \langle \boldsymbol{\phi}({\bf{x}}'),\boldsymbol{S}_{\omega}({\bf{x}},{\bf{x}}') \big \rangle,\boldsymbol{\phi}({\bf{x}}) \Big \rangle  \text{,}
\end{equation}
where $\boldsymbol{S}_{\omega}$, called the cross-spectral density tensor, is the Fourier transform pair of the statistically stationary space-time correlation tensor,

\begin{equation} \label{SPOD: CSD}
    \boldsymbol{S}_{\omega}({\bf{x}}_1,{\bf{x}}_2) = \int \limits_{-\infty}^{\infty}\boldsymbol{C}(\boldsymbol{x}_1,\boldsymbol{x}_2,\tau)e^{-i\omega \tau}d\tau 
\end{equation}
with
\begin{equation}
\boldsymbol{C}(\boldsymbol{x}_1,\boldsymbol{x}_2,\tau) = \mathbb{E}[\boldsymbol{q}(\boldsymbol{x}_1,t + \tau)\boldsymbol{q}^*(\boldsymbol{x}_2,t)] \text{.}
\end{equation}
Note that in statistically stationary flow, the correlation tensor $\boldsymbol{C}$ only depends on one time variable, $\tau$, which represents the difference in the two times in a more general $\boldsymbol{C}$ that arises in space-time POD.
\\

The mathematical structure of (\ref{SPODCorrelation}) is the same as in space-only POD, so to maximize (\ref{SPODCorrelation}) we solve the integral eigenvalue problem
\begin{equation}
    \int \limits_{\Omega}\boldsymbol{S}_{\omega}({\bf{x}}_1,{\bf{x}}_2)\boldsymbol{W}({\bf{x}}_2)\boldsymbol{\psi}_{\omega,k}({\bf{x}}_2)d{\bf{x}}_2 = \lambda_k\boldsymbol{\psi}_{\omega,k}({\bf{x}}_1) \text{,}
\end{equation}
where $\lambda_k$ is the $k^{\text{th}}$ largest eigenvalue of $\boldsymbol{S}_{\omega}$, and $\boldsymbol{\psi}_{\omega,k}$ is the $k^{\text{th}}$ SPOD mode at frequency $\omega$.

\subsection{With discrete data}
As we saw in \ref{susbsec:SOPOD_disc}, data can be used to compute approximate space-only POD modes by forming a discrete spatial correlation matrix ${\bf{C}}$. Analogously, approximate spectral POD modes can be computed by using data to form a discrete cross-spectral density matrix. To do this, instead of taking the Fourier transform of the space-time correlation matrix, the Wiener-Khinchin theorem is invoked, which states that the cross-spectral density tensor is equivalent to the correlation between points in the flow in Fourier space,
\begin{equation}
    \boldsymbol{S}_{\omega}(\boldsymbol{x}_1,\boldsymbol{x}_2) = \mathbb{E}[\hat{\boldsymbol{q}}_{\omega}(\boldsymbol{x}_1)\hat{\boldsymbol{q}}_{\omega}^*(\boldsymbol{x}_2)] \text{.}
\end{equation}
To find ${\bf{S}}_{\omega}$ in this way, many realizations of the flow are needed in the frequency domain, at a particular frequency: $\hat{\bf{Q}}_{\omega} = [\hat{\bf{q}}_{\omega,1},\hat{\bf{q}}_{\omega,2},\dots \hat{\bf{q}}_{\omega,m}]$, where $\hat{\bf{q}}_{\omega,i} \in \mathbb{C}^N$ is the spatially discretized representation of the Fourier transform of the $i^{th}$ flow realization, at frequency $\omega$. These frequency domain realizations come from applying the discrete Fourier transform to time series data, and to get more than one realiztion at each frequency from a single time series, the time series is broken up into $m$ (possibly overlapping) blocks, and the discrete Fourier transform is taken of each block \cite{Aaron18}.  With this data matrix of Fourier realizations, the cross-spectral density matrix can be approximated as
\begin{equation}
{\bf{S}}_{\omega} = \frac{1}{m}\hat{\bf{Q}}_{\omega}\hat{\bf{Q}}_{\omega}^{*}  \text{.}
\end{equation}
The product of the data matrix with its transpose approximates the discrete cross-spectral density tensor because it averages products of different components of $\hat{{\bf{q}}}_{\omega}$ over the realizations. The SPOD modes at frequency $\omega$ are the eigenvectors of this approximated cross-spectral density tensor multiplied by the weight matrix, 
\begin{equation}
{\bf{S}}_{\omega}{\bf{W}}{\bf{\Psi}}_{\omega} = {\bf{\Psi}}_{\omega}{\bf{\Lambda}}_{\omega}\text{,}
\end{equation}
where ${\bf{\Psi}}_{\omega} = [{\bf{\psi}}_{\omega,1}, {\bf{\psi}}_{\omega,2} \dots ]$, and ${\bf{\Lambda}}_{\omega}$ is the diagonal matrix of eigenvalues at frequency $\omega$. As before, the modes can be obtained using the SVD of the (weighted) data matrix of realizations in Fourier space,
\begin{equation}
    \frac{1}{\sqrt{m}}{\bf{W}}^{\frac{1}{2}}\hat{{\bf{Q}}}_{\omega} = {\bf{U}}_{\omega}{\bf{\Sigma}}_{\omega}{\bf{V}}^*_{\omega} \text{.}
\end{equation}
The modes and energies are then given by
\begin{equation}
    {\bf{\Psi}}_{\omega} = {\bf{W}}^{-\frac{1}{2}}{\bf{U}}_{\omega} \quad \text{and} \quad {\bf{\Lambda}}_{\omega} = {\bf{\Sigma}}^2_{\omega}\text{.}
\end{equation}
\\

\section{Space-time POD} \label{sec:STPOD}
In the previous two sections, we described space-only POD, where the goal is to find modes that optimally represent snapshots of the flow, and spectral POD, where the goal is to find modes that optimally represent the Fourier transform of the flow at each individual frequency. Next, we introduce a generalization of these methods called space-time POD, in which the goal is to find modes that optimally represent the flow evolution on a finite-time window $[0,T]$.
\\

\subsection{Formulation}
The goal is formalized in the same way as the previous two cases, by maximizing the expected value of the square of an inner product,

\begin{subequations} \label{STPOD: optimization}
\begin{equation}
    \lambda_{{\bf{xt}}} [\boldsymbol{\phi}({\bf{x}},t) ] = \frac{1}{\|\boldsymbol{\phi}({\bf{x}},t)\| }  \mathbb{E}\big[ \|\langle \boldsymbol{q}({\bf{x}},t), \boldsymbol{\phi} ({\bf{x}},t)\rangle_{\bf{xt}} \|^2 \big] \text{,}
\end{equation}
\begin{equation}
    \boldsymbol{\phi}_1({\bf{x}},t)  = \argmax \lambda_{\bf{xt}}[\boldsymbol{\phi}({\bf{x}},t)] \text{.}
\end{equation}
\end{subequations}
This time, however, our inner product acts over both space and time, 
\begin{equation} \label{STPOD: innerP}
    \langle \boldsymbol{q}_1({\bf{x}},t),\boldsymbol{q}_2({\bf{x}},t) \rangle_{\bf{xt}}  = \int \limits_{0}^{T}\int \limits_{\Omega}\boldsymbol{q}_2^{{\bf{*}}}({\bf{x}},t)\boldsymbol{W}({\bf{x}},t)\boldsymbol{q}_1({\bf{x}},t)d{\bf{x}}dt \text{.}
\end{equation}
\\
The expectation operator acts over finite time segments of the flow. These segments can be separate realizations of the flow or can be extracted from a single long time series if the flow is ergodic. The space-time optimization problem (\ref{STPOD: optimization}) has a different physical meaning than the space-only optimization and the spectral optimization problems (\ref{SOPOD:Optimization}) and (\ref{SPOD: optimization}), respectively. The desired mode here will be a function of both space and time, and its accuracy is measured by an inner product over both space and time as well, whereas, in the previous two cases, the modes were functions of only space. Mathematically, however, the problems are quite similar; instead of the spatial domain, $\Omega$, we now have the spatiotemporal domain $\mathcal{D} \equiv \Omega \times [0,T]$, and the modes are now a functions of both space and time, ${\bf{z}} \equiv [{\bf{x}},t]^T$. With these definitions, the space-time optimization problem can be written
\begin{subequations} 
\begin{equation} \label{Sop: maximization}
     \lambda_{\bf{xt}}[\boldsymbol{\phi}({\bf{z}})] =  \frac{1}{\| \boldsymbol{\phi}({\bf{z}})\| } \mathbb{E}\big[ \|\langle \boldsymbol{q}({\bf{z}}), \boldsymbol{\phi} ({\bf{z}})\rangle_{\bf{st}} \|^2 \big] \text{,}
\end{equation}
\begin{equation}
    \boldsymbol{\phi}_1({\bf{z}})  = \argmax \lambda_{\bf{xt}}[\boldsymbol{\phi}({\bf{z}})] \text{,}
\end{equation}
\end{subequations}
where the inner product is
\begin{equation}
      \langle \boldsymbol{q}_1({\bf{z}}),\boldsymbol{q}_2({\bf{z}}) \rangle_{\bf{xt}}  = \int \limits_{\mathcal{D}}\boldsymbol{q}_2^{{\bf{*}}}({\bf{z}})\boldsymbol{W}({\bf{z}})\boldsymbol{q}_1({\bf{z}})d{\bf{z}} \text{.}
\end{equation}
The maximization is exactly the same mathematically as (\ref{SOPOD:Optimization}), so we know the solution will be that the optimal space-time modes are the eigenfunctions of the correlation between two points in the spatiotemporal domain,
\begin{equation}
        \int \limits_{\mathcal{D}}\boldsymbol{C}({\bf{z}}_1,{\bf{z}}_2)\boldsymbol{W}({\bf{z}}_2)\boldsymbol{\phi}_{k}({\bf{z}}_2)d{\bf{z}}_2 = \lambda_k\boldsymbol{\phi}_{k}({\bf{z}}_1) \text{.}
\end{equation}
Rewriting this again in terms of space and time, the space-time POD modes satisfy
\begin{equation} \label{STPOD: cont eig}
    \int \limits_{0}^{T}\int \limits_{\Omega}\boldsymbol{C}({\bf{x}}_1,t_1,{\bf{x}}_2,t_2)\boldsymbol{W}({\bf{x}}_2,t_2)\boldsymbol{\phi}_{k}({\bf{x}}_2,t_2)d{\bf{x}}_2dt_2 = \lambda_k\boldsymbol{\phi}_{k}({\bf{x}}_1,t_1) \text{,}
\end{equation}
where the space-time correlation tensor is defined as 
\begin{equation}
    \boldsymbol{C}({\bf{x}}_1,t_1,{\bf{x}}_2,t_2) = \mathbb{E}[\boldsymbol{q}({\bf{x}}_1,t_1) \boldsymbol{q}({\bf{x}}_2,t_2)^{{\bf{*}}}] \text{.}
\end{equation}
\\

In section \ref{sec:Hankel}, we will show that the Hankel SVD modes constitute a discrete approximation of space-time POD modes. Before doing so, we first show that spectral and space-only POD can be recovered from space-time POD in the limits of long and short time windows, respectively. The latter result is novel to the best of our knowledge. 
\\

\subsection{Space-time POD on long times becomes spectral POD} \label{STPOD:Long time}
Here, we show that in the limit that $T \to \infty$, the space-time modes become spectral POD modes. This limit has been shown in the literature, e.g., in \cite{Lumley70,Aaron18}, but here we give a different proof. We restrict our attention to ergodic flows, though all that is needed for the results to be valid is that the flow is wide-sense statistically stationary, i.e., the first and second moments of the flow don't change with time. The integral eigenvalue problem that defines the modes when $T \to \infty$ is
\begin{equation} \label{eq:STPOD:inf limit}
    \int \limits_{-\infty}^{\infty}\int \limits_{\Omega}\boldsymbol{C}({\bf{x}}_1,t_1,{\bf{x}}_2,t_2)\boldsymbol{W}({\bf{x}}_2)\boldsymbol{\phi}({\bf{x}}_2,t_2)d{\bf{x}}_2dt_2 = \lambda\boldsymbol{\phi}({\bf{x}}_1,t_1) \text{.}
\end{equation}
For convenience, the bounds on time have been shifted to cover the entire real line. Because the flow is stationary, the correlation must only depend on the difference of the two times, not on the times themselves, so 
\begin{equation}
    \boldsymbol{C}({\bf{x}}_1,t_1,{\bf{x}}_2,t_2)  \to \boldsymbol{C}({\bf{x}}_1,{\bf{x}}_2,t_1-t_2)\text{.}
\end{equation}
Switching the order of integration, we have
\begin{equation} \label{eq:STPOD:}
    \int \limits_{\Omega}\int \limits_{-\infty}^{\infty}\boldsymbol{C}({\bf{x}}_1,{\bf{x}}_2,t_1-t_2)\boldsymbol{W}({\bf{x}}_2)\boldsymbol{\phi}({\bf{x}}_2,t_2)dt_2d{\bf{x}}_2 = \lambda\boldsymbol{\phi}({\bf{x}}_1,t_1) \text{,}
\end{equation}
which is a convolution of $\boldsymbol{C}$ with $\boldsymbol{\phi}_{k}$. Taking the Fourier transform of both sides and using the convolution theorem, we have
\begin{equation} \label{eq:STPOD:SPOD}
    \int \limits_{\Omega}\boldsymbol{S}_{\omega}({\bf{x}}_1,{\bf{x}}_2)\boldsymbol{W}({\bf{x}}_2)
    \hat{\boldsymbol{\phi}}_{\omega}({\bf{x}}_2)d{\bf{x}}_2 = \lambda\hat{\boldsymbol{\phi}}_{\omega}({\bf{x}}_1) \text{.}
\end{equation}
The modes must satisfy (\ref{eq:STPOD:SPOD}) for every $\omega$. SPOD modes, $\boldsymbol{\psi}_{\omega'}$ at some $\omega'$, are solutions as they are defined by this equation for $\omega = \omega'$, and are zero for $\omega \neq \omega'$. The SPOD mode with the $k^{\text{th}}$ greatest energy $\lambda$ (over the modes at all frequencies) is the $k^{\text{th}}$ space-time POD mode. 
\\

\subsection{Space-time POD on short times becomes space-only POD} \label{subsec:STPOD:short time}
Here, we show an analogous result for space-only POD. Specifically, we investigate what happens to the space-time POD modes in the limit that the time $T$ on which they evolve is short compared to the time scales of the flow. In this limit, the space-time correlation between two points $({\bf{x}}_1,t_1)$ and $({\bf{x}}_2,t_2)$ is only a function of their locations in space because $t_1$ and $t_2$ are negligibly different. Intuitively, we can now think of this integral over space and time as an integral over space multiplied by time because there is no time dependence in the integrand, so long as the weight doesn't depend on time (which would be an unusual choice). Rewriting condition (\ref{STPOD: cont eig}) with a constant weight and correlation, we have
\begin{equation}
        \int \limits_{0}^{T}\int \limits_{\Omega}\boldsymbol{C}({\bf{x}}_1,{\bf{x}}_2)\boldsymbol{W}({\bf{x}}_2)\boldsymbol{\phi}_{k}({\bf{x}}_2,t_2)d{\bf{x}}_2dt_2 = \lambda_k\boldsymbol{\phi}_{k}({\bf{x}}_1,t_1) \text{.}
\end{equation}
Because the left-hand side is not a function of $t_1$, the modes must indeed be constant in time, and the condition on the modes becomes an integral over space multiplied by time,
\begin{equation}
       T\int \limits_{\Omega}\boldsymbol{C}({\bf{x}}_1,{\bf{x}}_2)\boldsymbol{W}({\bf{x}}_2)\boldsymbol{\phi}_{k}({\bf{x}}_2)d{\bf{x}}_2 = \lambda_k\boldsymbol{\phi}_{k}({\bf{x}}_1) \text{.}
\end{equation}
This condition is exactly that of space-only POD, as shown in (\ref{SOPOD:integral_evp}). As a result, the space-time POD modes converge to the space-only POD modes and their energies are proportional to the space-only POD energies (with proportionality constant $T$).
\\

\subsection{Space-time POD with discrete data}
Given time series data on a spatial grid from an experiment or simulation, how do we compute the space-time POD modes? Above, we have seen that the solution to the space-time POD problem is given by eigendecomposition of the space-time correlation, which becomes a matrix when space and time are discrete. Building this matrix requires any ensemble of finite-time realizations of the flow. If the time $T$ that the space-time POD modes are to evolve on is represented by $d$ time steps with our temporal discretization, then one realization of the flow can be written as a vector with $d$ snapshots stacked on top of one another,

\begin{equation}
{\bf{y}}  = \begin{bmatrix} 
    {\bf{q}}_1\\{\bf{q}}_2\\ 
    \vdots\\
   {\bf{q}}_d
    \end{bmatrix} \in \mathbb{R}^{Nd} \text{.}
\end{equation}

With many of these realizations, we can approximate the correlation matrix ${\bf{C}} \in \mathbb{R}^{Nd\times Nd}$ in the same manner as before,
\begin{equation} \label{eqn:GDPOD_CorrelationTEns}
{\bf{C}} = \frac{1}{m}{\bf{Y}}{\bf{Y}}^{{\bf{*}}} \quad \text{with} \quad {\bf{Y}}  = [{\bf{y}}_1,{\bf{y}}_2,\dots {\bf{y}}_m] \in \mathbb{R}^{(Nd) \times m}\text{.}
\end{equation}
\\
Analogous to the space-only POD and SPOD cases, this approximation works because the correlation between two points in the spatio-temporal domain is approximated as the product of the flow at those points averaged over the realizations. We will see later in section \ref{sec:exploiting_ergodicity}, however, that if the flow is assumed to be ergodic and the realizations are all to be formed from one long time series, this construction of the correlation matrix does not fully exploit the ergodicty of the system, and it is possible to create a more accurate approximation of the correlation. The discretized POD modes are obtained as
\begin{equation}
    {\bf{C}}{\bf{W}}{\bf{\Phi}} = {\bf{\Phi}}{\bf{\Lambda}}\text{,} 
\end{equation}
where ${\bf{W}} \in \mathbb{R}^{Nd \times Nd}$ is the discretized weight matrix. 
\\

Similar to the space-only POD and SPOD cases, the modes on this domain can also be obtained from the SVD of the weighted data matrix,
\begin{equation} \label{eq:D:Weight_SVD} 
    \frac{1}{\sqrt{m}}{\bf{W}}^{\frac{1}{2}}{\bf{Y}} = {\bf{U}}{\bf{\Sigma}}{\bf{V}}^* \text{.}
\end{equation}
The space-time POD modes and energies are then
\begin{equation} \label{eq:D:Weight}
{\bf{\Phi}} = {\bf{W}}^{-\frac{1}{2}}{\bf{U}} \quad \text{and} \quad {\bf{\Lambda}} = {\bf{\Sigma}}^2 \text{.}
\end{equation}
\\

\section{Hankel singular vectors are space-time POD modes} \label{sec:Hankel}
Here, we demonstrate the connection between space-time POD modes and Hankel modes by showing that one way of computing the unweighted space-time POD modes is to take the SVD of a Hankel matrix. As we will show, this connection provides insight into the sense in which the Hankel SVD modes are optimal and the impacts of the time step, rows, and columns on the properties of the modes.
\\

For ergodic systems, the flow realizations needed to construct the approximate correlation matrix can be extracted from a single long time series. Given a time series of snapshots, ${\bf{q}}_1,  {\bf{q}}_2, \dots , {\bf{q}}_d, \dots, {\bf{q}}_{m+d-1}$, we may extract $m$ realizations of length $d$ by creating the first realization ${\bf{y}}_1 = [{\bf{q}}_1^T  {\bf{q}}_2^T, \dots, {\bf{q}}_d^T]^T$, then advancing one time step over to create the second realization, ${\bf{y}}_2 = [{\bf{q}}_2^T  {\bf{q}}_3^T, \dots, {\bf{q}}_{d+1}^T]^T$, and so on. Stacking these as columns in the data matrix, the result is a block Hankel matrix,
\begin{equation} \label{HankelFlows}
    {\bf{H}} = [{\bf{y}}_1,{\bf{y}}_2, \dots {\bf{y}}_m] = 
    \begin{bmatrix} {\bf{q}}_1 & {\bf{q}}_2 & \dots & {\bf{q}}_m \\ {\bf{q}}_2 & {\bf{q}}_3  & \dots & {\bf{q}}_{m+1}\\ 
    \vdots & \vdots & \ddots & \vdots\\
    {\bf{q}}_d &  {\bf{q}}_{d+1}  & \dots  & {\bf{q}}_{m+d-1} 
    \end{bmatrix}
    \text{.}
\end{equation}
Comparing (\ref{HankelFlows}) and (\ref{eqn:GDPOD_CorrelationTEns}) we see that the Hankel matrix ${\bf{H}}$ provides one way of generating the data matrix ${\bf{Y}}$. Therefore, using (\ref{eq:D:Weight_SVD}) and (\ref{eq:D:Weight}) and taking the SVD of the (weighted) Hankel matrix,
\begin{equation} \label{eq:Hankel_SVD_Weight2} 
    \frac{1}{\sqrt{m}}{\bf{W}}^{\frac{1}{2}}{\bf{H}} = {\bf{U}}{\bf{\Sigma}}{\bf{V}}^* \text{,}
\end{equation}
we obtain the space-time POD modes as
\begin{equation} \label{eq:Hankel_SVD_Weight1}
{\bf{\Phi}} = {\bf{W}}^{-\frac{1}{2}}{\bf{U}} \quad \text{and} \quad {\bf{\Lambda}} = {\bf{\Sigma}}^2 \text{.}
\end{equation}
If the weight ${\bf{W}}$ is uniform, the modes come directly from the Hankel SVD,
\begin{equation}
    {\bf{\Phi}} = {\bf{U}} = \text{lsv}({\bf{H}}) \text{.}
\end{equation}
That is, the left singular vectors of the Hankel matrix give a discrete approximation of the space-time POD modes \textit{in the case of a uniform weight}. 
\\

This connection to space-time POD implies the sense in which the Hankel modes are optimal: they minimize the expected value of the square norm of the projection error, where the norm is the $\ell_2$ norm, i.e., the weight is uniform. Uniformity is a reasonable choice for the weight matrix if the spatial grid is uniform. However, if the grid is non-uniform, a uniform weight matrix corresponds to a non-uniform continuous weight. For example, if the flow data is defined on a cylindrical grid, a uniform weight matrix will bias accuracy toward the center of the domain over the outside. If the desired weight is non-uniform, as, e.g., is likely to be the case if the grid is non-uniform in space, then the discrete space-time POD modes can be obtained using equations (\ref{eq:Hankel_SVD_Weight1}) and (\ref{eq:Hankel_SVD_Weight2}).
\\

The connection we've drawn between space-time POD and Hankel singular modes provides insight into the distinct meaning of rows and columns of the Hankel matrix, which we investigate next.
\\

\subsection{Height of the Hankel matrix corresponds to $T$} \label{subsec:Height of Hankel}
\begin{figure} [ht]
    \centering
    \includegraphics[scale = 1]{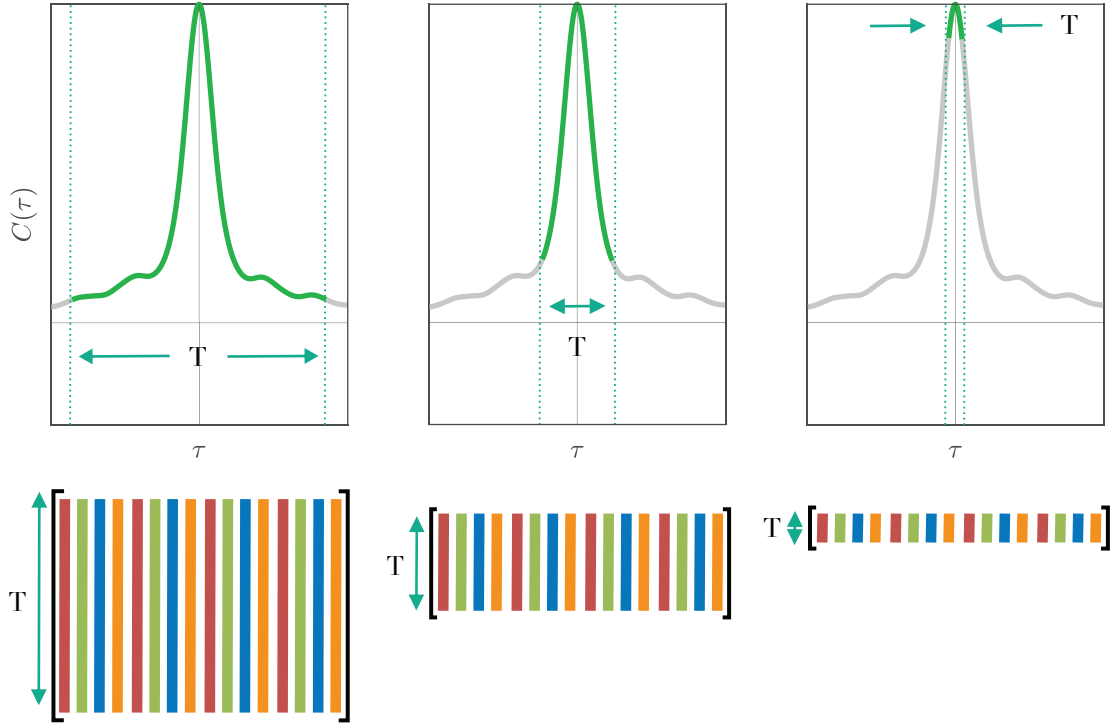}
    \caption{The height of the Hankel matrix $d$, along with the time step $\Delta t$ between snapshots, determines the time window $T$ over which the space-time correlation $C$ of the system is accounted for in the Hankel singular modes. Shorter Hankel matrices, therefore, severely truncate the correlations, whereas taller Hankel matrices retain more of these correlations.}
    \label{fig:Correlation Truncation}
\end{figure}

The connection between space-time POD and singular modes of the Hankel matrix reveals two interpretations of the number of rows $d$ in the Hankel matrix. First, it determines, for a fixed time step $\Delta t$, the time window $T$ on which the space-time POD modes optimally represent the flow. The window also depends on the time step, specifically, $T= (d-1) \Delta t$. Second, the significance of the choice of $d$ can also be understood in terms of the space-time correlation matrix implied by the Hankel matrix --- the Hankel singular modes will account for a section of the space-time correlation of the system of width $T$, as shown in Figure \ref{fig:Correlation Truncation}. The properties of the correlation outside of this window are discarded; they play no role in the definition of the modes and cannot be represented by the modes. 
\\

Two interesting limits of this truncation occur when $(d-1)*\Delta t << \tau_C$ and when $(d-1)\Delta t >>\tau_C$, where $\tau_C$ is the maximum correlation time of the flow between any two locations in space. When the height of the Hankel matrix corresponds to a time much smaller than the correlation time, we should expect to get spatial modes identical to those from discrete space-only POD for the same reason discussed in section~\ref{subsec:STPOD:short time}: the discrete space-time POD modes are eigenvectors of the correlation matrix, and if $(d-1)*\Delta t = T << \tau_C$, then this correlation matrix is constant in time. In this limit, the correlation between two points in space is independent of the time difference, so 
\begin{equation}
    {\bf{C}}_{ij} = {\bf{C}}_{i+aN,j+bN} \text{,}
\end{equation}
where $N$ is the number of spatial gridpoints and $a$ and $b$ are any integers, so long as the indices are valid. Incrementing one of the indices by $N$ corresponds to moving one time step forward but not shifting in space, but we have assumed that the correlation matrix is constant in time. It can therefore be written in block form as 
\begin{equation}
    {\bf{C}} = \begin{bmatrix}
    {\bf{C}}^{\bf{x}} && {\bf{C}}^{\bf{x}}&& \dots && {\bf{C}}^{\bf{x}} \\
    {\bf{C}}^{\bf{x}} && \ddots && && \vdots\\
    \vdots \\
    {\bf{C}}^{\bf{x}} && \dots && && {\bf{C}}^{\bf{x}}
    \end{bmatrix} \text{,}
\end{equation}
where ${\bf{C}}^{\bf{x}} \in \mathbb{R}^{N \times N}$ is the discrete correlation matrix at no time lag, as defined in \ref{eqn:SOPOD_CorrelationTEns}. ${\bf{C}}$ has only $N$ non-zero eigenvalues, and these eigenvalues are $d\Lambda^{\bf{x}}$, the space-only eigenvalues scaled by the number of time steps. The associated eigenvectors of ${\bf{C}}$ are
\begin{equation}
    {\bf{\phi}}_k  = \begin{bmatrix}
    {\bf{\phi}}_k^{\bf{x}}\\
    {\bf{\phi}}_k^{\bf{x}}\\
    \vdots\\
    {\bf{\phi}}_k^{\bf{x}}
    \end{bmatrix} \text{,}
\end{equation}
where ${\bf{\phi}}_k^{\bf{x}}$ is the $k^{\text{th}}$ eigenvector of the space-only correlation ${\bf{C}}^{\bf{x}}$. Since each $N\times 1$ block of ${\bf{\phi}}_k$ is the same, the space-time POD modes are constant in time and have the spatial form of space-only POD. Thus, we've recovered space-only POD in the short-time limit as we might expect from the discussion of the continuous case in section \ref{subsec:STPOD:short time}. 
\\

Above, as well as in section \ref{subsec:STPOD:short time}, we assumed no temporal variation of the correlation function. Gibson et al. \cite{Gibson92} showed that for the special case of a scalar Hankel matrix, if small variations in time to the correlation matrix are allowed, then the resulting POD basis is the Legendre polynomials. In intuitive terms, the Legendre polynomials appear because, over short times, the flow will be well approximated by its Taylor series truncated at a few terms, and each term is much more important than the next. Therefore, with $k$ modes, the basis should span the set of degree $k-1$ polynomials and also be orthogonal with respect to the inner product. If the inner product is uniform, then the basis satisfying these conditions is the Legendre polynomials. On a vector time series, we observe a full set of modes with no time dependence, consistent with the order zero Legendre polynomial, followed by modes with negligible energy whose evolution follows subsequent Legendre polynomials. That is, the spatial variance is much greater, and therefore more important to capture, compared to the temporal variation. This confirms that our analysis assuming no temporal variation, and the resulting conclusions that the space-time modes converge to space-only modes, hold. 
\\
\\
The other limit, $(d-1)\Delta t = T >>\tau_C$, occurs when the height of the Hankel matrix corresponds to a time much longer than any correlation time in the flow. From section \ref{STPOD:Long time}, we would expect that the modes oscillate with a pure frequency as SPOD modes do. Indeed, taking the SVD of a Hankel matrix becomes a discrete Fourier transform as the height of the Hankel matrix grows to infinity \cite{Broomhead86,Bozzo10}. This limit has been used in the context of time-delay embedding \cite{Kaiser20}.
\\

Finally, we note that the time step between rows (along columns) signifies the temporal discretization of the continuous eigenvalue problem (\ref{STPOD: cont eig}) which defines the modes. Next, we leverage the distinct meaning of the time step between rows and columns to approximate the modes at reduced cost. 
\subsection{Width of the Hankel matrix informs convergence} \label{subsec:Columns}
The connection with space-time POD also reveals the impact of the width of the Hankel matrix on its singular modes. The convergence of the discrete space-time POD modes is determined by the accuracy of the approximation of the correlation matrix formed from the data. Each element of the correlation is proportional to the sample average of products of one entry in the column with another, averaged over the columns of the Hankel matrix, 
\begin{equation}
    {\bf{C}}_{ij} = ({\bf{H}}{\bf{H}}^*)_{ij}  = \sum_{k = 1}^{m}{\bf{H}}_{ik}{\bf{H}}_{kj} \text{.}
\end{equation}
Thus, the width of the Hankel matrix $m$ determines the number of realizations that contribute to approximating the correlation.
\\

However, because two adjacent columns in the Hankel matrix are only one time step apart, their contributions to each element of the correlation matrix are far from independent, so the accuracy of the correlation is not just a function of the number of columns. Indeed, one could imagine a Hankel matrix that has many columns but only represents data over a short time during which the underlying system does not fully explore its phase space. There are two criteria necessary for the modes to be accurate: there must be enough columns in the Hankel matrix so that their sample average converges, and the data must be representative of the underlying attractor. The latter may require the data to be taken over a long time, leading to a very wide Hankel matrix, thus a costly SVD. We point out here that the time step between columns signifies the time between two successive realizations, in contrast to the time step between rows. 
\\

\begin{figure}
\captionsetup[subfloat]{labelformat=empty}
\subfloat[]{
    \begin{overpic}[trim={-0.4cm 0 0 0},clip]{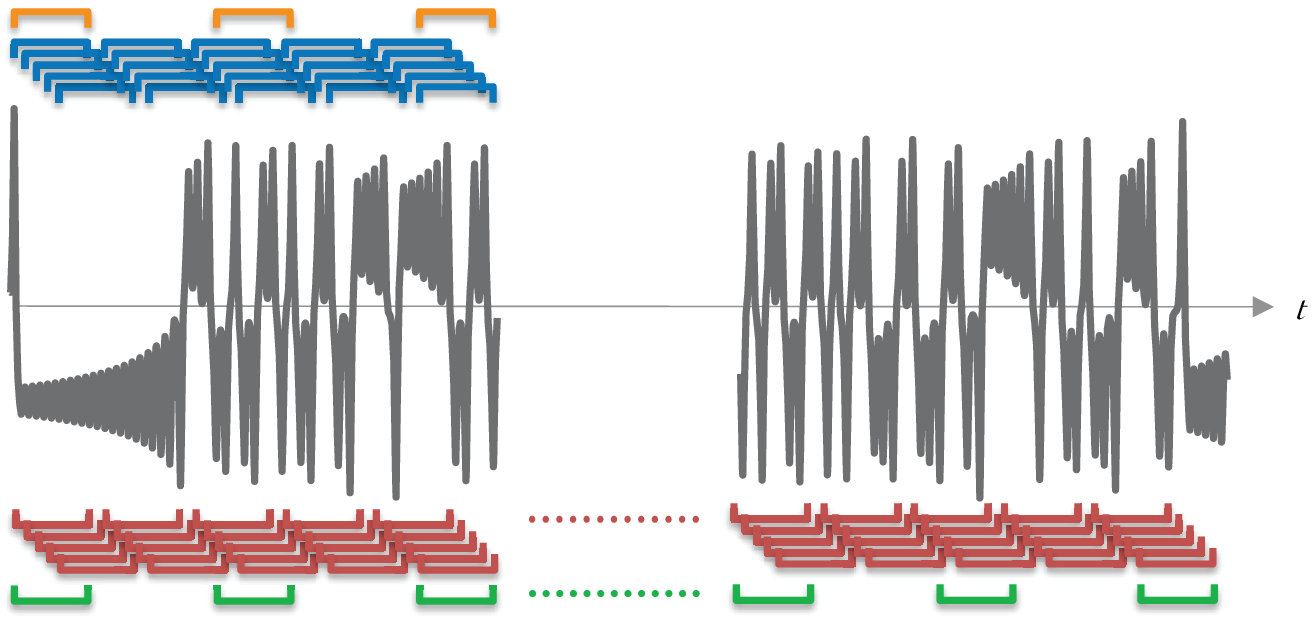}
    \put(-2,44){(a)}
    \end{overpic}
    }\\
    \centering
    \input{Correlation_Graphic}
    \includegraphics{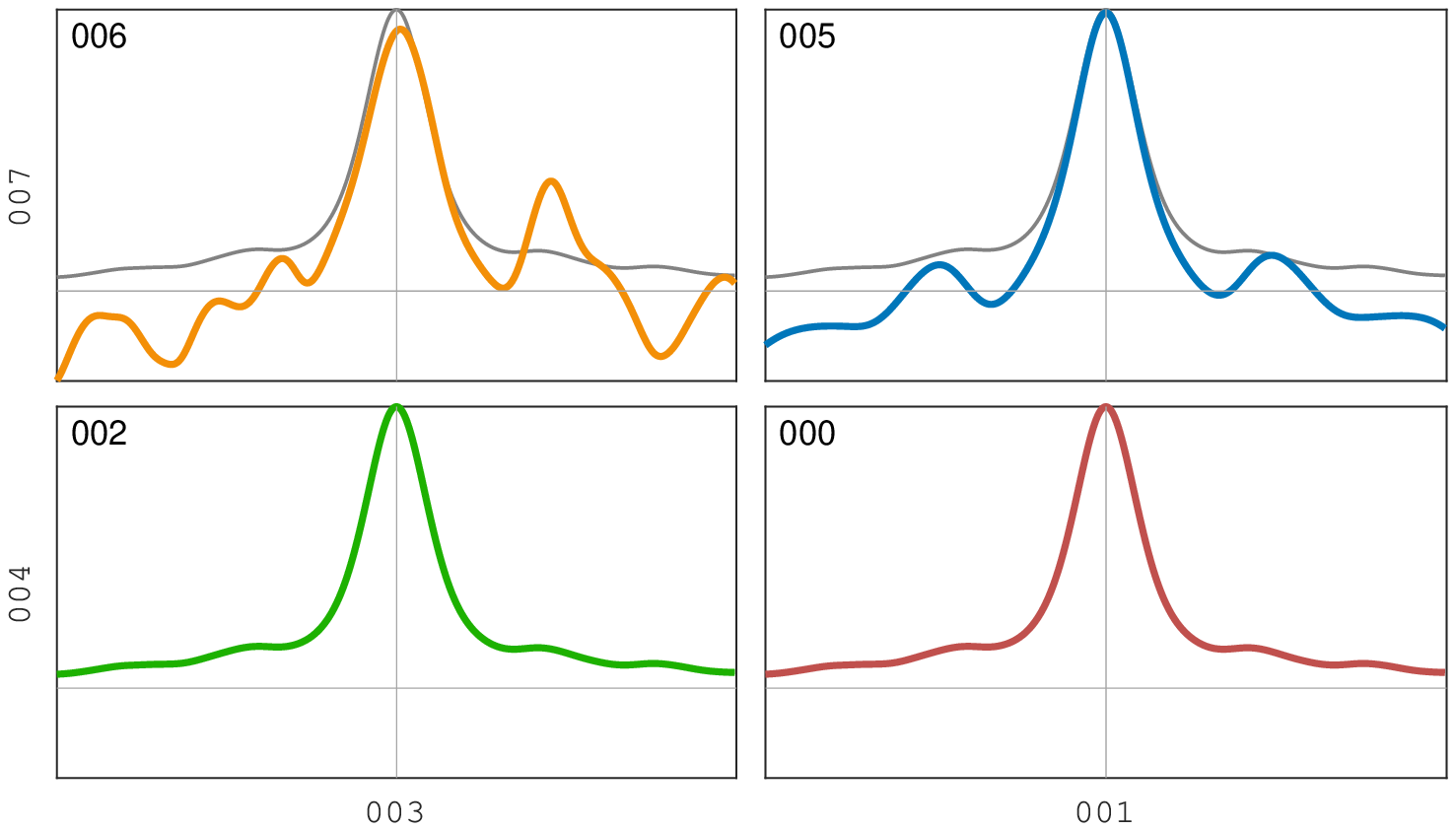}
    \caption{Uncorrelated columns (left) and Hankel (right) approaches to approximating the correlation matrix from short (b,c) and long (d,e) time series. (a) A time series used to sample temporal realizations with different samplings. (b) Correlations from uncorrelated columns using a short time series. (c) Correlations from the Hankel matrix using the short times series. (d) Correlations from uncorrelated columns using a long time series. (e) Correlations from a Hankel matrix using the long time series. Though both (d) and (e) are accurate, the correlations in (d) come at significantly lower computational cost.}
    \label{correlation cartoons}
\end{figure}
Alternatively, we may throw out many of the columns of the Hankel matrix, so that the columns are less correlated, then take the SVD of this matrix,

\begin{equation} \label{eq:Q_UQ}
{\bf{Q}}_{UC} =  \begin{bmatrix} {\bf{q}}_1 & {\bf{q}}_{s+1} & {\bf{q}}_{2s+1} & \dots & {\bf{q}}_m \\ {\bf{q}}_2 & {\bf{q}}_{s + 2}  & {\bf{q}}_{2s+2} & \dots & {\bf{q}}_{m+1}\\ 
    \vdots & \vdots & \vdots & \ddots & \vdots\\
    {\bf{q}}_d &  {\bf{q}}_{s + d}  & {\bf{q}}_{2s+d} & \dots  & {\bf{q}}_{m+d-1} 
    \end{bmatrix}
    \text{.}
\end{equation}
Here, $s$ represents the spacing in time between columns, and the Hankel matrix is recovered if $s = 1$. ${\bf{Q}}_{UC}{\bf{Q}}_{UC}^*$ still forms an approximation of the correlation matrix, so its singular vectors will approximate the space-time POD modes. If ${\bf{H}}$ is formed from data representative of the attractor, then ${\bf{Q}}_{UC}$ will be as well, and if we retain enough columns, the correlation matrix will be accurate. An appropriate choice of $s$, one which negligibly impacts the accuracy of the modes, depends on the amount of data available, the system dimension, the embedding dimension, and the time step. The Hankel matrix is the data matrix with the most columns possible from the data, and in practice it is likely overkill, and overly costly, if the time series is long. If the data is lacking, however, the Hankel matrix will produce a more accurate correlation and hence more accurate modes than any other spacing. 
\\

This strategy of removing most of the columns provides a substantial cost reduction. The time complexity of the SVD scales quadratically with the smaller dimension of the matrix and linearly with the larger one. Therefore, this strategy reduces the computation time by a factor of $s^2$ if the Hankel matrix is taller than it is wide, which is usually the case for fluids and other high dimensional problems, and by a factor of $s$ otherwise. Critically, this reduction is achieved without changing the window $T$ or the time step $\Delta t$ along the columns, such that the meaning of the modes discussed in section~\ref{subsec:Height of Hankel} and the temporal discretization of the integral in (\ref{STPOD: cont eig}) remain unchanged. 
\\

A cartoon of the correlations produced from uncorrelated columns vs. Hankel columns with different amounts of data is shown in Figure~\ref{correlation cartoons}. The top two correlation graphics show how the uncorrelated columns (left) and Hankel (right) approaches might approximate the correlation with limited data. In this case, the Hankel matrix is likely to produce more accurate correlations than the matrix with uncorrelated columns because the latter may have too few to converge averages. However, the limited data is not representative of the underlying attractor, so although sampling this data more and more (Hankel approach) and computing the correlations may converge, it will not converge to the true correlations. When the time series is long enough to be representative of the attractor (bottom), both approaches will produce accurate correlations so long as enough columns are included in the uncorrelated matrix to converge averages. We also note that in the case of a short time series, we derive a method in section~\ref{sec:exploiting_ergodicity} which produces a more accurate approximation of the correlations than the Hankel approach with the same data.
\\

\subsection{Summary of the connection}
\begin{figure}
    \centering
    \begin{overpic}[trim={-0.4cm 0 0 -0.4cm},clip]{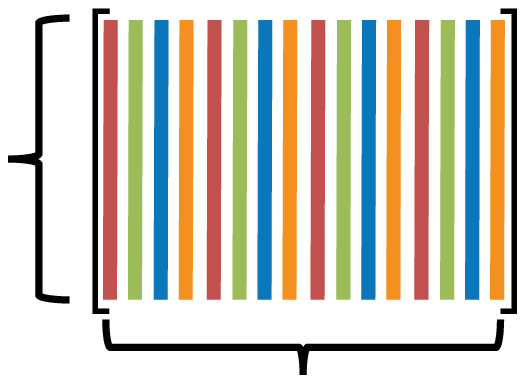}
    \put(8,23){$
        \bigintss \limits_{0}^{T} \bigintss \limits_{\Omega}\boldsymbol{C}({\bf{x}}_1,t_1,{\bf{x}}_2,t_2)\boldsymbol{W}({\bf{x}}_2,t_2)\boldsymbol{\phi}_{k}({\bf{x}}_2,t_2)d{\bf{x}}_2dt_2 = \lambda_k\boldsymbol{\phi}_{k}({\bf{x}}_1,t_1)
    $}

    \linethickness{2pt}
    \put(15,57){\color{myred} \vector(-1.9,-10){5.1}}
    \put(38,40.5){\color{myred} \vector(-1.4,-1){21}}
    \put(72,27.5){\color{myred} \vector(0,1){28}}
    \put(29.8,73){Hankel Matrix}
    \put(-10,61){Height corresponds}
    \put(-10,57){to $T$}
    \put(60,58){\huge{$= {\bf U\Sigma V}^*$}}
    \put(38.4,38.5){Width corresponds to}
    \put(38.6,35.5){accuracy of correlation}
    \put(73,44){${\bf u}_k \to \boldsymbol{\phi}_k$}
    \put(73,39){$\sigma_k \to \sqrt{m\lambda_k}$}

    \end{overpic}
   \vspace*{-25mm}
    \caption{The dimensions of the Hankel matrix can be interpreted in terms of the integral space-time eigenvalue problem (\ref{STPOD: cont eig}). The height of the Hankel matrix corresponds to $T$, the length on which the space-time POD modes optimally represent the flow. The width, though there is subtlety about independence, informs how accurately the correlation matrix is approximated, and hence the accuracy of the modes.}
    \label{fig:Hankel Schematic}
\end{figure}

To summarize, we have shown that Hankel singular modes constitute an approximation of space-time POD modes. The block Hankel matrix is formed as a data matrix whose columns are discretized temporal flow realizations, and when multiplied by its conjugate transpose, gives the correlation matrix from space-time POD. The eigenvectors of this matrix, which are the space-time POD modes, are the left singular vectors of the Hankel matrix due to the well-known equivalence of the SVD and eigendecomposition. This connection sheds light on the sense in which the Hankel modes are optimal: they are optimal in the case of a uniform weight, which is often an undesirable norm. Understanding Hankel singular modes as approximate space-time POD modes also reveals that the time step between rows and columns of the Hankel matrix need not be the same. The time step along a column represents the discretization of the integral eigenvalue problem \ref{STPOD: cont eig}, which defines the continuous modes. The time step along a row corresponds to the time between successive flow realizations. This distinction implies that while a short time step along columns is desirable, a short time step along rows may not be, so we leverage this distinction to reduce the computation time by removing many of the columns. In many practical cases, this negligibly impacts the accuracy of the modes but significantly reduces the time for the SVD. The time window corresponding to the height of the Hankel matrix $T = (d-1)\Delta t$ is the time over which the Hankel modes optimally represent the flow as well as the window of space-time correlations accounted for in calculating the modes. Finally, in the limits of short and tall Hankel matrices, corresponding to evolution on short and long time windows, the Hankel modes become discrete space-only POD modes and discrete SPOD modes, respectively. Much of the discussion in this section is summarized in Figure \ref{fig:Hankel Schematic}.
\\

\section{Fully exploiting ergodicity: more accurate modes} \label{sec:exploiting_ergodicity}
As discussed in section \ref{sec:Hankel}, the Hankel matrix implies a particular approximation of the space-time correlation, which in turn determines the approximation of space-time POD modes provided by singular modes of the Hankel matrix. In this section, we show how a more accurate approximation of the correlation matrix can be constructed by fully exploiting the ergodicity of the underlying system, which produces more accurate modes if there is a shortage of data, as we will demonstrate later in section \ref{subsec:Exploiting_Ergodicity}. This approach has been used in the context of SSA \cite{VAUTARD89}, but to our knowledge, has never been employed for vector-valued data.
\\

Ergodicity is assumed in constructing the Hankel matrix: instead of taking columns from different sample paths of the flow, columns are formed from different sections in one long sample path, so there is an assumption that these two are equivalent. However, in constructing the Hankel matrix (or a down-sampling thereof) and then taking its SVD, the assumed ergodicity is not fully exploited. To see this, look at the correlation ${\bf{C}}  = {\bf{H}}{\bf{H}}^* \in \mathbb{R}^{Nd \times Nd}$. Broken up into its spatial block structure, this correlation is written
\begin{equation}
    {\bf{C}} = \begin{bmatrix}
    {\bf{C}}_{00} && {\bf{C}}_{01}&&  {\bf{C}}_{02}&& \dots && {\bf{C}}_{0d} \\
   {\bf{C}}_{01}^T && {\bf{C}}_{11}&&  {\bf{C}}_{12}&& \dots && {\bf{C}}_{1d}\\
   {\bf{C}}_{02}^T && {\bf{C}}_{12}^T&&  {\bf{C}}_{22}&& \dots && {\bf{C}}_{2d}\\
    \vdots && \vdots && \vdots && \ddots && \vdots \\
   {\bf{C}}_{0d}^T && {\bf{C}}_{1d}^T&&  {\bf{C}}_{2d}^T&& \dots && {\bf{C}}_{dd}\\
    \end{bmatrix} \text{,}
\end{equation}
where the matrix
\begin{equation}
    {\bf{C}}_{ij} \in \mathbb{R}^{N \times N}
\end{equation}
represents the correlations between all points in the spatial domain between times $i$ and $j$. Because the flow is ergodic, these correlations should only depend on the difference in times, e.g., ${\bf{C}}_{12}$ should be equal to ${\bf{C}}_{56}$, but by calculating this correlation with the Hankel matrix, this will not normally be the case. Indeed, writing each of these correlation blocks in terms of elements of the time series,
\begin{equation} \label{eq:Ci_Hankel}
    {\bf{C}}_{ij} = \sum_{k = 1}^{m}{\bf{q}}_{i+k-1}{\bf{q}}_{j + k - 1}^* \text{,}
\end{equation}
we see, for example, that the ($1,1$) block of ${\bf{C}}$ does not depend on the last element of the time series while the ($d,d$) block does, so they will not be equivalent as we know they should be from ergodicity. Instead, the fully converged correlation should have the symmetric block Toeplitz structure,
\begin{equation} \label{eq:Toeplitz Correlation}
    \tilde{{\bf{C}}} = \begin{bmatrix}
    {\bf{C}}_{0} && {\bf{C}}_{1}&&  {\bf{C}}_{2}&& \dots && {\bf{C}}_{d} \\
   {\bf{C}}_{1}^T && {\bf{C}}_{0}&&  {\bf{C}}_{1}&& \dots && {\bf{C}}_{d-1}\\
   {\bf{C}}_{2}^T && {\bf{C}}_{1}^T&&  {\bf{C}}_{0}&& \dots && {\bf{C}}_{d-2}\\
    \vdots && \vdots && \vdots && \ddots && \vdots \\
   {\bf{C}}_{d}^T && {\bf{C}}_{d-1}^T&&  {\bf{C}}_{d-2}^T&& \dots && {\bf{C}}_{0}\\
    \end{bmatrix} \text{.}
\end{equation}
That the (block) diagonals of ${\bf{C}} = {\bf{H}}{\bf{H}}^*$ are not constant indicates that ergodicity is not fully exploited in constructing the correlation.
\\

To calculate each correlation ${\bf{C}}_i$ from a time series, ${\bf{q}}_1, {\bf{q}}_2, \dots {\bf{q}}_{m+d-1}$, we simply take the sample average of the product of all points in the time series $i$ time steps apart,
\begin{equation} \label{eq:Ci_Toeplitz}
    {\bf{C}}_i = \frac{1}{m+d-i}\sum_{k = 1}^{m+d-i} {\bf{q}}_k{\bf{q}}^*_{k+i} \text{.}
\end{equation}
After building $\tilde{{\bf{C}}}$ as in (\ref{eq:Toeplitz Correlation}), the modes are obtained as $\tilde{{\bf{C}}}{\bf W} {\bf \Phi }= {\bf{\Phi}}{\bf{\Lambda}}$. The difference in accuracy of the modes can be substantial when data is in short supply and when m/d is not large, as we demonstrate in section~\ref{subsec:Exploiting_Ergodicity}.
\\

The entries of the Toeplitz matrix (\ref{eq:Ci_Toeplitz}) are formed using more terms than are used in forming the Hankel-based correlations (\ref{eq:Ci_Hankel}). The ratio depends on the entry calculated but is proportional to $\frac{m}{m+d}$; therefore as $m/d$ increases, the advantage in the Toeplitz approach diminishes.
\\

The relative cost of computing modes using the Hankel and Toeplitz approaches depends on the relative size of $N d$ and $m$. The time scaling for the Hankel approach comes from the SVD in both cases, which is $\mathcal{O}(N^2d^2m)$ for $Nd<m$ and $\mathcal{O}(Ndm^2)$ otherwise. The time scaling for the Toeplitz approach comes from passing through the time series to calculate all blocks of the correlation and from taking the eigendecomposition of the correlation. The former scales as $\mathcal{O}((m+d)dN^2)$ and the latter as $\mathcal{O}(N^3d^3)$. The ratio of the scalings of the Toeplitz and Hankel algorithms is thus
\begin{equation}
    \begin{cases}
    \max(\frac{1}{d},\frac{Nd}{m}) \quad \text{if $Nd<m$,} \\ \\
   \big(\frac{Nd}{m}\big)^2 \quad \text{if $Nd>m$.}
    \end{cases}
\end{equation}

In fluids applications, $Nd$ is usually much greater than $m$, so this approach will scale much worse than the Hankel approach. However, in many other dynamics applications, including many classic applications of the Hankel matrix, $N$ is small, often $1$, and $d$ substantially smaller than $m$ \cite{Brunton16,Dylewsky22}, so this approach both leads to more accurate modes and better scaling. In short, the Toeplitz method scales more favorably than building the Hankel matrix and taking the SVD when $Nd<m$, but when $Nd>m$ the algorithm is slower. 
\\

\section{Numerical experiments: Lid-driven cavity flow} \label{sec:example}
In this section, we demonstrate our theoretical results with data from a 2D lid-driven cavity flow at Reynolds number $Re = 22{,}000$ \cite{Cazemier98}. First, we compare the convergence of modes with a true Hankel matrix against a data matrix whose columns are uncorrelated, showing that with the same number of columns (flow realizations) the uncorrelated matrix generates better modes. We then show a related, but more practical, example. Given an $m_{H}$-column Hankel matrix, one can form an uncorrelated data matrix by sampling $m$ of its columns. We show that for $m$ significantly smaller than $m_H$ the modes are nearly identical, which can be used to reduce computational cost. Next, we demonstrate that space-time POD converges to space-only POD on very short time intervals, and that the modes become Fourier in time on very long time intervals, indicating convergence to spectral POD. Finally, we show that with limited data, exploiting the ergodicity of the system produces more accurate modes than the Hankel approach.
\\

\subsection{Simulation description} \label{sec:Numerical Experiments_Simulation Description}
\begin{figure}
    \centering
    \begin{overpic}[scale = 0.55]{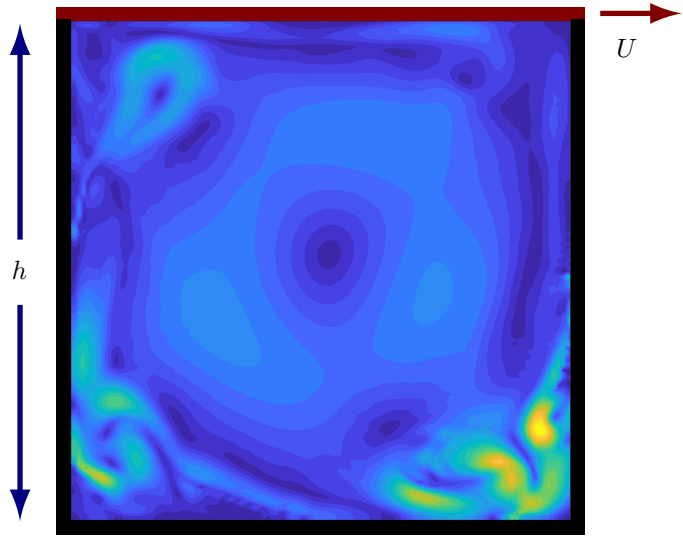}
      %\put(20,69.4){\color{myred}\rule{0.56\textwidth}{5pt}}
      \put(19.5,69.4){\color{myred}\rule[0pt]{199.6pt}{5pt}}
      \put(19.5,8.1){\color{black}\rule[0pt]{5pt}{190pt}}
      \put(19.5,6.5){\color{black}\rule[0pt]{199.6pt}{5pt}}
      \put(82.45,8.1){\color{black}\rule[0pt]{5pt}{190pt}}
      
      \linethickness{2pt}
      \put(86,70.3){\color{myred} \vector(1,0){10}}
      %\put(8,70.3){\color{myred} \vector(1,0){10}}
      \put(15,42.6){\color{myblue} \vector(0,1){26.5}}
      \put(15,34.6){\color{myblue} \vector(0,-1){26.5}}
      \put(14,37.8){$h$}

      \put(88,65){$U$}
    \end{overpic}
    \caption{Schematic of the lid-driven cavity flow. The lid moves to the right and drives the flow. Contours show a snapshot of the turbulent kinetic energy.}
    \label{fig:snapshot}
\end{figure}

We generate data for square lid-driven cavity flow at $Re = \frac{Uh}{\nu} = 22{,}000$, where $U$ is the speed of the lid, $h$ is the height of the square cavity, and $\nu$ is the viscosity of the fluid. We solve the incompressible Navier-Stokes equations using a Crank-Nicolson method for the viscous term and an Adams-Bashforth method for the nonlinear term \cite{michio20}. The domain is discretized with $N_x = N_y = 120$ grid-points. Data is generated by starting the simulation with zero velocity everywhere except the top, running until the initial transients have vanished and the statistics become stationary, and then collecting the time series data. The time is nondimensionalized so that in one unit of time the lid slides the cavity width, and the simulation time step is $5 \times 10^{-4}$. We sample this data with $\Delta t = 0.1$. Figure \ref{fig:snapshot} shows a snapshot of the flow to illustrate the setup and spatial scale of variation of the flow.
\\

\subsection{Spaced columns yield better modes}
In section \ref{subsec:Columns}, we argued that an $m$-column data matrix with uncorrelated columns will produce a more accurate correlation matrix, and hence more accurate modes, than an $m$-column Hankel matrix, when measured against converged modes. Here, this is demonstrated with the lid-driven cavity flow: we form two data matrices --- one Hankel, and the other uncorrelated, with the structure in (\ref{eq:Q_UQ}). We choose a downsampled version of $u$, the $x$ component of velocity, as our observable of the state ${\bf{q}}$, with $24$ points in each direction. We choose $T = 2$, so that the modes evolve on the convective time scale of the flow, and $d =\frac{T}{\Delta t} + 1 = 21$. Each column of the matrices is an element of $\mathbb{R}^{Nd}$, where $Nd = 24 \times 24 \times 21 = 12{,}096$. The space-time POD modes are obtained as the left singular vectors of the two matrices, and their convergence is evaluated by taking the square of their inner product with a fully converged mode, which is obtained from a data matrix with $40{,}000$ columns, all far apart in time. The variation in different realizations of the chaotic system leads to variation in the accuracy of the modes calculated from them. To account for this, we repeat the process of finding the modes $400$ times for each $m$ and take the mean square inner product with the fully converged mode.
\\

\begin{figure}
    \centering
    \input{Uncorrelated_vs_Correlated}
    \includegraphics[scale = 1]{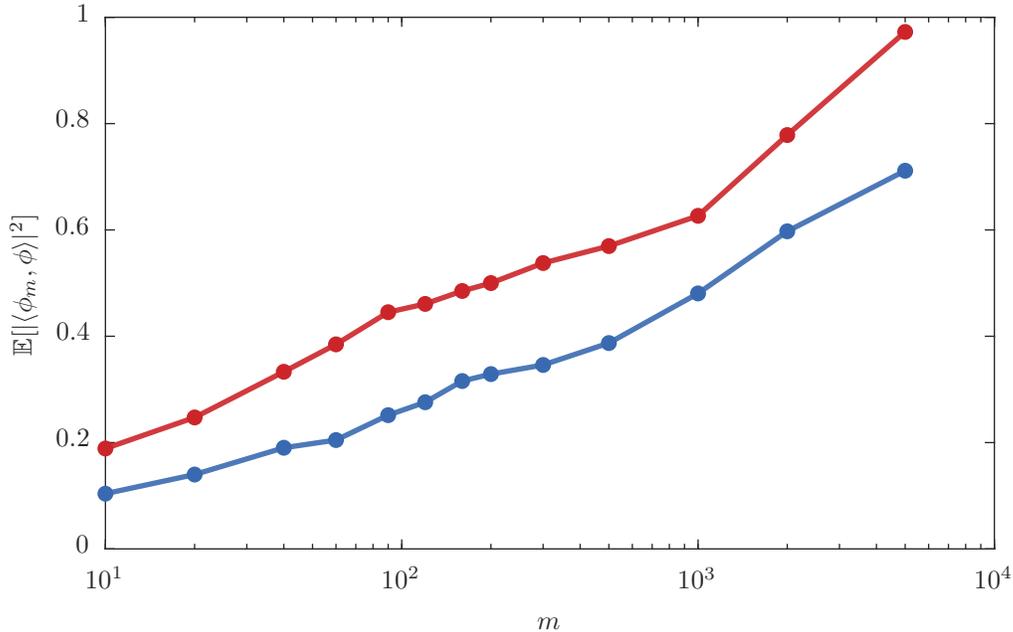}
    \caption{Convergence of the modes with the width $m$ of the matrix: (blue) Hankel matrix; (red) uncorrelated data matrix. The modes from the uncorrelated data matrix converge faster than those from the Hankel matrix. This improved convergence is observed because the columns of the uncorrelated data matrix represent data from more of the attractor than the columns of the Hankel matrix of the same size.}
    \label{fig:Uncorrelated vs Hankel}
\end{figure}
Figure \ref{fig:Uncorrelated vs Hankel} shows the results. As expected, the modes from the uncorrelated data matrix outperform those from the Hankel matrix. The $x$-axis is logarithmic, and the two curves are roughly a constant horizontal distance from one another. This means that to achieve the same accuracy as an uncorrelated data matrix with $m$ columns, a Hankel matrix must have $cm$ columns, where $c$ is some constant independent of $m$; here it is roughly $8$, though this number depends sensitively on the parameters used.
\\

\subsection{Faster computation by sampling the Hankel matrix}
Choosing between taking the SVD of a Hankel matrix or an uncorrelated data matrix of the same size may be unrealistic --- simulation or experimental data can be in short supply, and it is not trivial to simply generate the extra data needed to construct the uncorrelated data matrix. A more relevant question might be: given a time series of snapshots, should one form the entire Hankel matrix or can many of the columns be omitted? The latter has the advantage of reducing computational cost, in terms of both CPU time and memory requirements, while the former may be more accurate. The accuracy loss depends on how many columns are removed, but we demonstrate here that the loss of accuracy can be negligible for orders of magnitude computational speedup. 
\\

Again, we choose $T = 2$, so $d = 21$. Starting with a Hankel matrix with $m_H = 5000$ columns, we form a reduced data matrix by retaining only $m$ columns with the maximum possible spacing. For example, if $m = 51$, we use every $100^{\text{th}}$ column of the Hankel matrix to form the data matrix. Removing columns changes the time step along rows but not along columns, and our analysis in section \ref{subsec:Columns} indicates that this only changes the rate of convergence of the modes, not what the modes converge to. We take the SVD of both the Hankel and data matrices and compare each mode from the data matrix to its respective Hankel matrix mode using the square of the inner product; values of $1$ and $0$ indicate identical modes and orthogonal modes, respectively. We repeat this process $200$ times to see how well, on average, the modes from the much smaller data matrix approximate the modes from the full Hankel matrix. 
\\

\begin{figure}
    \centering
    \input{moreDataEffect}
    \includegraphics[scale = 1]{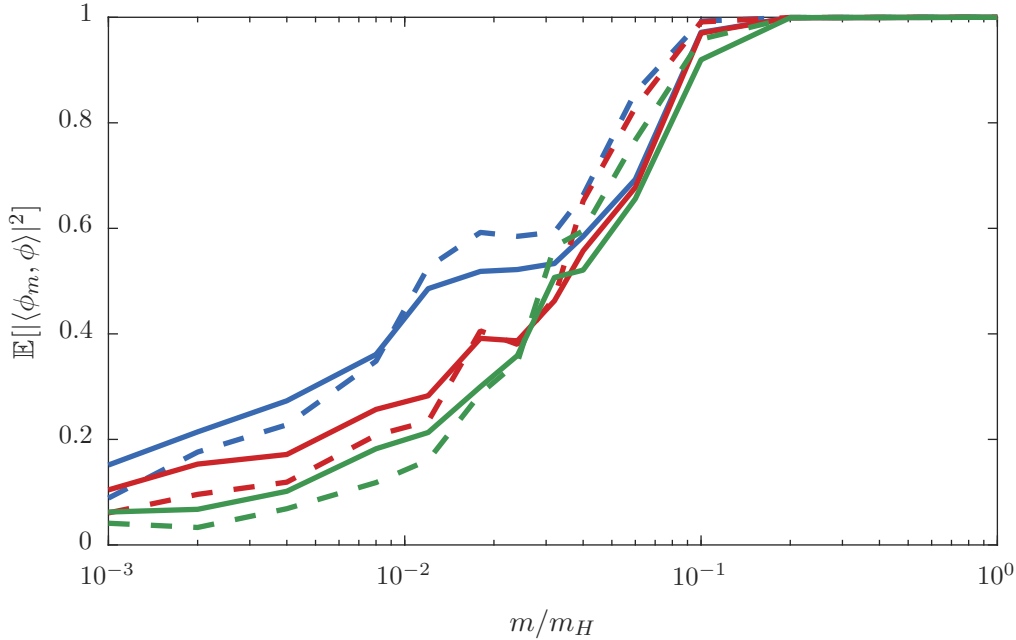}
    \caption{Similarity between modes from a downsampled Hankel matrix, retaining $m/m_H$ columns, and the Hankel modes. Both the mean (solid) and the median (dashed) are reported for the $200$ trials. The first mode (blue) and the second mode (red) are more similar than the third mode (green). For the parameters used, retaining $1/10$ and $1/5$ of the columns yields modes that capture over $97\%$ and $99.9\%$ of the energy of the Hankel modes, respectively, for both the first and second modes, on average (mean).}
    \label{fig:More Data Effect}
\end{figure}

Figure \ref{fig:More Data Effect} shows this accuracy metric as a function of $m$ for the first three modes from each matrix, reporting both the mean (solid) and the median (dashed) of the $200$ trials. Both the first and second modes from the $m = 500$ data matrix capture over $97\%$ of the energy of the corresponding Hankel modes (mean), and over $99\%$ (median). For $m = 1000$, the first and second modes both capture over $99.9\%$ of the corresponding Hankel modes. Even these modest reductions in $m$ offer significant computational savings, scaling here like the square of the reduction in $m$.
\\

The ratio of the number of columns in the sampled Hankel matrix to that of the full Hankel matrix needed to achieve a given accuracy is not constant. For example, if the time step $\Delta t$ is smaller, a sparser sampling of the Hankel matrix can be used to achieve the same accuracy because $k$ columns of the Hankel matrix represent less of the underlying flow. The speedup gained by sampling the Hankel matrix can be significant; for problems with a complicated spatial domain, the height of the data matrices, $Nd$, will likely be much larger than the number of columns, $m$. Therefore the time complexity of the SVD is $m^2$, and reductions in $m$ save considerable time. We reiterate that this cost savings is enabled by the insight, provided by the connection between Hankel singular modes and space-time POD, that the time step along rows and columns have different meanings and need not be equal.

\subsection{Convergence to space-only POD for short times}
\begin{figure}[!hbt]
    \centering
    \input{ST_vs_POD}
    \includegraphics[scale = 1]{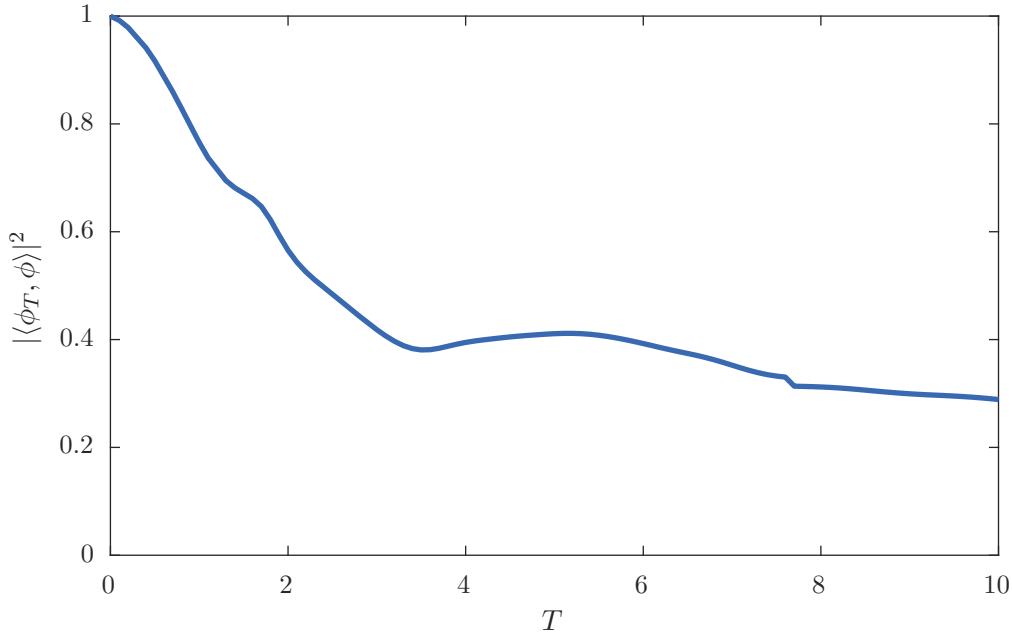}
    \caption{Convergence of space-time POD modes $\phi_{T}$ to space-only POD modes $\phi$ for small $T$. The spatial dependence of the first space-time mode at each $T$ is compared to the first space-only mode using the space-only inner product, averaged over the time evolution of the space-time mode.}
    \label{fig:ST_POD_to_POD:plot}
\end{figure}

Here, we demonstrate the convergence of the space-time POD modes to space-only POD modes as the time interval they are defined on approaches zero, $T \to 0$. Space-only POD modes are defined on $\Omega$ and space-time POD modes are defined on $\Omega \times [0,T]$, so in order to compare the two, we take the square of the space-only inner product, averaged over time. In Figure \ref{fig:ST_POD_to_POD:plot}, we see that the spatial part of the space-time modes indeed converges to space-only modes as $T \to 0$. In fact, the square inner product of the two remains near unity even when the space-time modes are defined over a time on the order of the convective time-scale; when $T=1$ it is $0.937$, and it only drops below $0.9$ when $T>1.4$.
\\

\subsection{Convergence to spectral POD for long times}
\begin{figure}[!hbt]
    \centering
    \input{Spectral_convergence}
    \includegraphics[scale = 1]{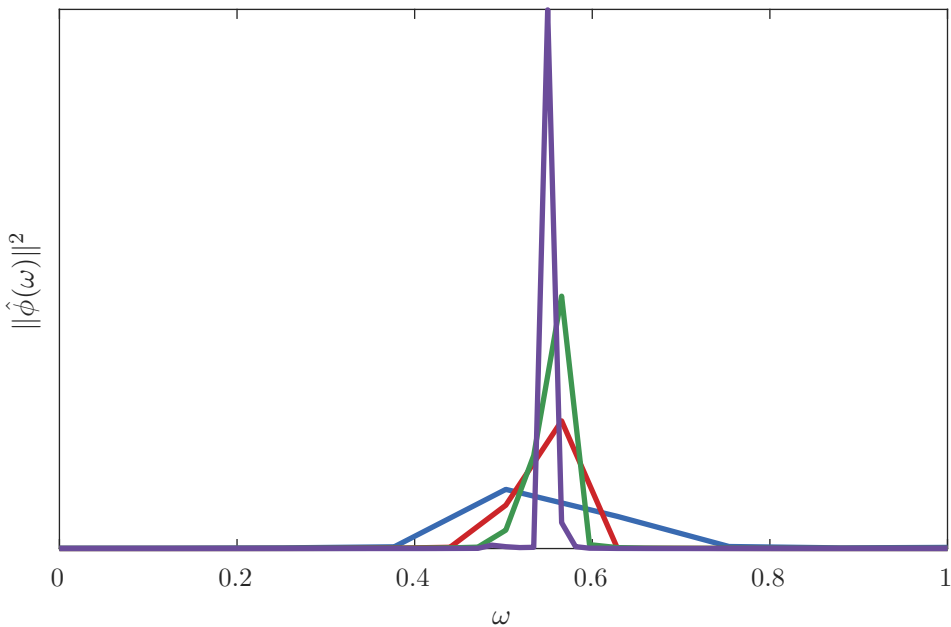}
    \caption{Convergence of space-time POD to spectral POD, visualized in terms of the power spectral density of the space-time POD modes for various values of $T$: (blue) $T = 50$; (red) $T = 100$; (green) $T = 200$; (purple) $T = 400$. As $T$ is increased, more and more of the energy of the modes is contained to an increasingly narrow band of frequencies, indicating convergence toward a discrete frequency, consistent with the time dependence of SPOD modes.}
    \label{fig:STPOD to SPOD}
\end{figure}

Here, we demonstrate that as the time interval becomes long, the time dependence of the modes becomes Fourier. Modes are computed via the SVD of a data matrix with large $T$ and $800$ well-spaced columns. We do this for different values of $T$ and examine the frequency content of the modes. The expectation is that the modes will become delta functions in the frequency domain as $T$ increases since we know that SPOD modes, which have purely Fourier time dependence, are the limits of space-time POD for long time windows. The spectral content of the space-time modes may vary over the spatial extent of the mode, so we use the norm of their Fourier component,
\begin{equation}
        \|\hat{\phi}(\omega)\|^2 = \langle \hat{\phi}(\omega,{\bf{x}}),\hat{\phi}(\omega,{\bf{x}})\rangle \text{,}
\end{equation}
which is equivalent to the power spectral density of the mode integrated over the domain, as a measure of their overall spectral content. Figure~\ref{fig:STPOD to SPOD} demonstrates the convergence with increasing $T$ of the spectra of the leading mode to a delta function in frequency.

\subsection{Fully Exploiting Ergodicity}\label{subsec:Exploiting_Ergodicity}
Finally, we demonstrate that the modes from the Toeplitz correlation matrix offer a significant improvement over those from the Hankel data matrix. We extract the modes by taking the SVD and eigendecomposition of the Hankel and Toeplitz matrices, respectively, and calculate the energy of the flow captured by the first few modes of both methods. We repeat this many times to generate PDFs of the energy captured, observing that the Toeplitz modes capture more energy. We also compare the modes from both methods to fully converged Hankel modes and tabulate the results for different parameters, again finding that the Toeplitz modes are more accurate.
\\

Energy of the space-time POD modes, (\ref{STPOD: optimization}), is an appropriate metric for evaluating the performance of the two methods for obtaining modes because it is the quantity optimized in the definition of the POD problem. Analogous to (\ref{sec:SOPOD:eq:Energy_corr_eig}), we may rewrite the energy of a mode in terms of its projection coefficients onto the exact modes,
\begin{equation}
    \lambda[{\bf{\phi}}] = \sum_{k = 1}^{\min (m,Nd)}c_k^2 \lambda_k \quad \text{with} \quad c_k = \langle \phi,\phi_k \rangle / \|\phi\| \text{,}
\end{equation}
where ${\bf{\phi}}_k$ is the $k^{\text{th}}$ space-time POD mode. We calculate the energy of the modes produced by the Hankel and Toeplitz methods using the above formula, where the `true' modes are calculated using a Hankel matrix built from a $40\text{,}000$-long time series with the same $N$ and $d$.
\\

Given $m$ and $d$, and an $m+d-1$-long time series, we compute the modes using the $Nd\times m$ Hankel matrix and the $Nd \times Nd$ Toeplitz matrix and calculate their energies. Of course, the accuracy of the modes from either method depends on the time series used, so we repeat this process $1000$ times, generating PDFs for different choices of $m$, $d$, and $N$. A PDF with significant support only near unity indicates a high probability of producing an accurate mode.
\\

\begin{figure}[!hbt]
    \input{PDFs_6plts}
    \centering
    \input{PDFs_6plts}
    \includegraphics[scale = 1]{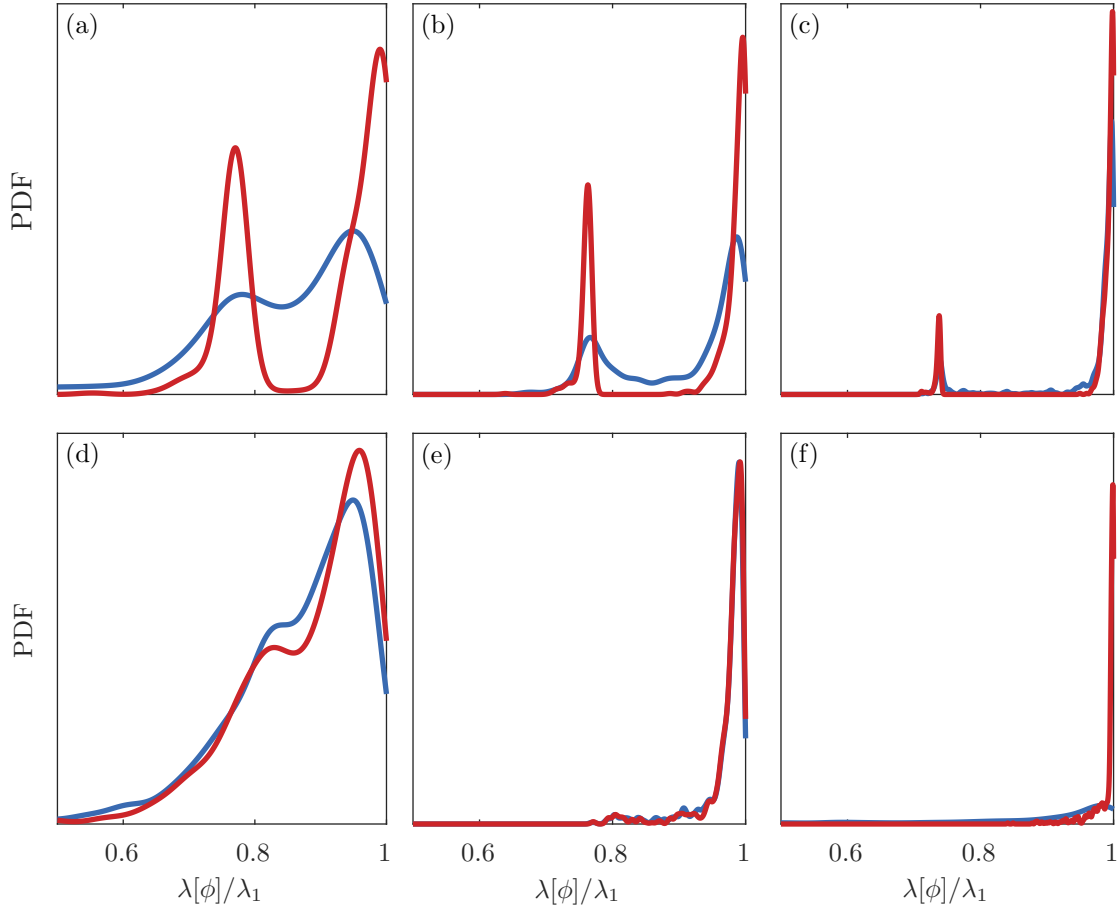}
    \caption{PDFs of energy captured by the first modes of the Hankel (blue) and Toeplitz (red) methods for various choices of $m$, $d$, and $N$. Both methods use the same time series of length $m+d-1$ to approximate the modes. The PDFs are calculated with $1000$ samples and the energies are divided by the maximum possible energy of one mode. These parameter values are, (a): $m = 30$, $d = 30$, $N = 1$. (b): $m = 100$, $d = 30$, $N = 1$. (c): $m = 300$, $d = 30$, $N = 1$. (d): $m = 300$, $d = 30$, $N = 3$. (e) $m = 600$, $d = 30$, $N = 3$. (f) $m = 10$, $d = 10$, $N = 1$.}
    \label{fig:6plts}
\end{figure}
Figure~\ref{fig:6plts} shows the PDF for the energy captured by the first mode from each method as a fraction of the energy of the true optimal mode ($\lambda_1$). Beginning with (a) in Figure \ref{fig:6plts}, we see that the mean of the modes from the Toeplitz method is higher than that of the modes from the Hankel method. Both methods give a bimodal distribution, and these peaks represent instances where the method mistakes the true second mode for the first mode. There is a deep trough between the peaks of the Toeplitz method, indicating that it produces a mode close to the first or to the second true mode rather than some combination thereof. The Hankel method displays some of this behavior, but it is significantly less pronounced, indicating that it produces something between the two true modes more often than the Toeplitz method. This behavior is also present in (b), (c), and (d). 
\\

As expected, for a fixed value of $d$ and $N$, increasing $m$ yields a higher probability of more accurate modes, and the modes from both methods in (c) are more accurate than those in (b), which are more accurate than those in (a). We also see that as $\frac{m}{d}$ increases, the difference in the accuracy of the two methods decreases, and the PDFs of the two methods become closer going from (a) to (b) and from (b) to (c). We see this behavior because for $\frac{m}{d} >> 1$, all of the elements on the same diagonal of ${\bf{C}} = {\bf{H H}}^*$ are sums of mostly the same products of terms from the time series, so ${\bf{C}}$ is nearly Toeplitz, whereas for $\frac{m}{d} \approx 1$ this is not true. In (e), $\frac{m}{d} = 50$, and the difference between the methods is hard to see by looking at the PDFs, though there is still some difference in the means and medians, especially for the second mode, as reported in Table~\ref{Table:Medians}.
\\

For fixed values of $d$ and $m$, increasing $N$ worsens the accuracy of the modes, and the modes in (d) are less accurate on average than those in (b). Similarly, increasing $d$ worsens the accuracy of the modes, as can be seen by comparing (a) and (f). In fact, the modes in (f) are more accurate from both methods than those in (a) despite being calculated with lower $m$. Increasing $N$ or $d$ decreases the quality of the modes because it causes them to be in a higher dimensional space, which makes them more difficult to approximate. If $m$ is large enough, however, (signifying a long time series from which to calculate the modes), both methods will converge for any $N$ and $d$. 
\\

\begin{figure}[ht]
    \centering
    \input{PDFs_multiplemodes}
    \includegraphics[scale = 1]{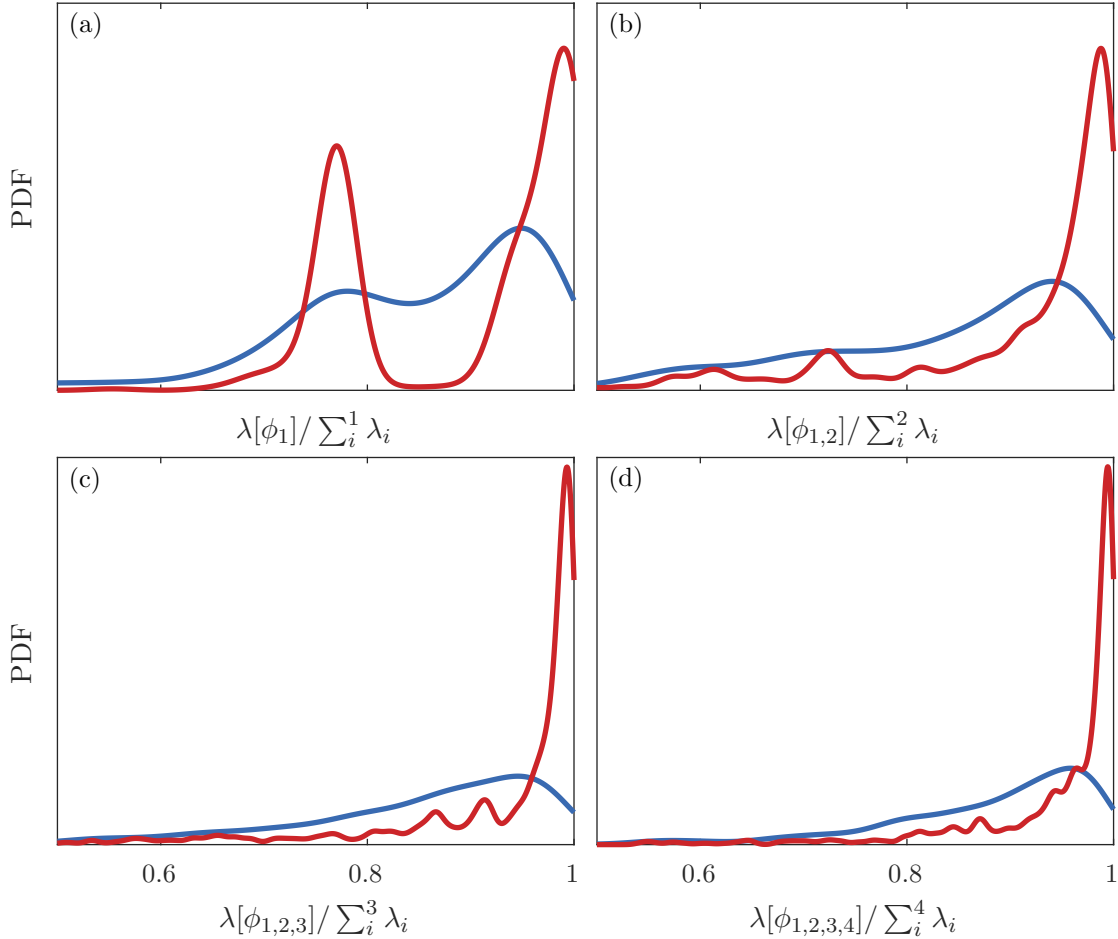}
    \caption{PDFs of the energy captured by the first four modes of the Toeplitz and Hankel methods, normalized by the maximum possible energy by four modes. These are calculated for $m = 30$, $d = 30$, $N = 1$. The difference is more severe for the latter modes, which we observed for the majority of the parameter choices that we tested.}
    \label{fig:multiple_modes}
\end{figure}
To evaluate the accuracy of the latter modes, say the $k^{\text{th}}$, we add the energies of the modes up to the $k^{\text{th}}$. We do this instead of simply calculating $\lambda[\phi_k]$ because this would erroneously reward switching the first and the $k^{\text{th}}$ mode. We denote the sum of the energies as
\begin{equation}
    \lambda[{\bf{\phi}}_{1,2,\dots , k}] = \sum_{i = 1}^k\lambda[{\bf{\phi}}_i] \text{,}
\end{equation}
and this quantity is interpreted as the energy captured by the first $k$ modes. Again, we generate PDFs of these energies using $1000$ time series of length $m+d-1$, computing modes using both methods, and calculating the energy they capture. Figure \ref{fig:multiple_modes} shows the results for the first four modes using $m = 30$, $d = 30$, $N = 1$. The energy captured by the first two Toeplitz modes is significantly greater than that of the first two Hankel modes. The same can be said about the first three and first four Toeplitz modes. These observations hold for the majority of parameter combinations.
\\

A second metric for evaluating the accuracy of the modes is to compare them directly to fully converged modes, via the square inner product. Again a value near unity indicates an accurate mode, though a value near zero for, e.g., the third mode does not indicate that that mode captures no energy, just that it captures energy orthogonal to that of the converged third mode. We calculate these converged modes from the $40\text{,}000$-long time series. For each choice of $m$, $d$, and $N$, we calculate the modes from the Hankel and Toeplitz matrices, take the inner product with the converged modes, and repeat $1000$ times. We record the mean and median of these inner products for the first three modes in Table \ref{Table:Medians}.
\\
 
\begin{table}[!hbt]
\centering
\begin{tabular}{ |p{3cm}|p{1.6cm}||P{2.8cm}|P{2.8cm}|P{2.8cm}|  }
 \hline
 Parameters& Method & \makecell{Mode 1\\ median / mean} & \makecell{Mode 2\\ median / mean} & \makecell{Mode 3\\  median / mean}\\
 \hline
 $m = 30$, $d = 30$, & Hankel   & $0.603$ / $0.532$    & $0.308$ / $ 0.375$ &   $0.223 $ / $0.320$ \\\cline{2-5}
 $N = 1$ & Toeplitz &   $0.930$ / $0.614 $ & $0.877 $ / $0.538$ & $0.808$ / $0.604$\\
 \hline
 \hline
 $m = 100$, $d = 30$, & Hankel   & $0.912$ / $0.672 $     & $0.885 $ / $0.651$ &   $0.928$ / $0.845$\\\cline{2-5}
 $N = 1$ & Toeplitz &   $0.978$ / $0.696  $  & $0.986 $ / $0.690$ & $0.957$ / $0.886$\\
 \hline
 \hline
 $m = 300$, $d = 30$, & Hankel   & $0.984 $ / $0.861$     & $0.989$ / $0.863$ &   $0.983$ / $0.970$\\\cline{2-5}
 $N = 1$ & Toeplitz &   $0.992$ / $0.880$  & $0.997$ / $0.884 $ & $0.986$ / $0.974$\\
 \hline
 \hline
 $m = 300$, $d = 30$, & Hankel   & $0.697$ / $0.571  $     & $0.637$ / $0.540$ &   $0.551$ / $0.508$\\\cline{2-5}
 $N = 3$ & Toeplitz &   $0.787$ / $0.601 $  & $0.732$ / $0.584$ & $0.591$ / $0.539$\\
 \hline
 \hline
 $m = 600$, $d = 30$, & Hankel   & $0.964$ / $0.891  $     & $0.954$ / $0.883$ &   $0.874$ / $0.754$\\\cline{2-5}
 $N = 3$ & Toeplitz &   $0.966$ / $0.904 $  & $0.956$ / $0.897 $ & $0.879$ / $0.755$\\
 \hline
 \hline
 $m = 10$, $d = 10$, & Hankel   & $0.919$ / $0.792$    & $0.678$ / $0.608 $&   $0.639$ / $0.567$\\\cline{2-5}
 $N = 1$ & Toeplitz &   $0.992$ / $0.880$  & $0.997$ / $0.884$ & $0.986$ / $0.974$\\
 \hline
\end{tabular}
\caption{Median and mean accuracy of Hankel modes vs. Toeplitz modes calculated from $1000$ time series of the most energetic $N$ points of the lid-driven cavity flow. Accuracy is measured as the average square inner product against a `true mode.' When the data is severely limited ($m$ is relatively small), the difference is more pronounced.}
\label{Table:Medians}
\end{table}

As $d$ increases, the number of realizations $m$ needed to get accurate modes also increases for both methods. Similarly as $N$ increases, $m$ also must increase to retain accuracy. Intuitively, this tells us that we need more data to find a mode in a higher dimensional space. The improvement of the Toeplitz method decays with $m/d$, because if this quantity is large, the difference between ${\bf{H}}{\bf{H}}^*$ and the Toeplitz $\tilde{\bf{C}}$ becomes small, as discussed in section~\ref{sec:exploiting_ergodicity}. Again, we note that the Toeplitz modes are significantly more accurate when $m$ is not significantly greater than $d$. Also, there is not a big difference when $m$ is very large, because here the modes are converged. Because a large $N$ necessitates a large $m$ for accuracy, $m/d$ will be large and thus the improvement gained by the Toeplitz method is marginal for large $N$. However, these are the cases for which the time complexity of the Toeplitz method is already prohibitive.
\\
 
\section{Conclusions} \label{sec:conclusions}
We have demonstrated that the singular modes of the Hankel matrix are a discrete approximation of the continuous space-time POD modes. That is, in the language of Dylewsky et al. \cite{Dylewsky22}, principal component trajectories are in fact an approximation of classical space-time POD modes defined on a finite temporal window. This connection is made by recognizing that the Hankel matrix multiplied by its conjugate transpose produces an approximation of a space-time correlation matrix, and the eigendecomposition of this matrix gives discrete space-time POD modes. We are able to gain useful insights into the Hankel modes by analyzing this correlation matrix rather than the Hankel matrix itself. This analysis makes clear the sense in which the modes are optimal: they are optimal if the weight matrix is uniform, which is often not desirable due to a variety of reasons including grid non-uniformity, the need to weight certain flow variables differently, or preference for certain regions in space. We derive a formula for the case of a non-uniform weight. The analysis of the correlation matrix also allows us to distinguish between the meaning of the time step along rows vs. along columns: the time step along columns corresponds to the temporal discretization of the space-time integral eigenvalue problem (\ref{STPOD: cont eig}), while the time step along rows corresponds to the time between successive flow realizations used to approximate the correlation matrix. Viewed in this way, it is clear that these time steps need not be the same. We leverage this insight to ease computation, omitting many of the columns of the Hankel matrix then taking the SVD of this new data matrix with uncorrelated columns. This produces modes with the same meaning at lower cost. We also show that two popular versions of POD, space-only POD and spectral POD, are recovered in the limits that the columns of the Hankel matrix represent short and long times, respectively. Finally, we improve the Hankel modes by fully exploiting the assumption of ergodicity to form a more accurate correlation matrix with the same data. This correlation matrix has Toeplitz structure owing to the fact that temporal correlations should only depend on the time separation, not on the times explicitly. We show that, especially if data is limited, the modes obtained as the eigenvectors of this Toeplitz matrix can be substantially more accurate than the Hankel modes. 
\\

We envision that this work will make Hankel-SVD-based methods applicable to new problems. In particular, the strategy of omitting many of the columns of the Hankel matrix may enable application of these methods to high dimensional systems such as those found in fluid dynamics. Likewise, the improved modes obtained by exploiting ergodicity may aid applications with severely limited time series data, such as expensive simulations. Furthermore, the theoretical connection we have established between space-time POD and Hankel SVD could lead to additional algorithmic improvements that further broaden the problems to which these methods are applied. Finally, given its connection to widely used Hankel methods, we hope this work will increase interest in space-time POD. 
\\

A number of questions remain. First, how does one, a priori, choose the number of columns to retain in the data matrix (how many to throw out from the Hankel matrix)? This depends on the data available, $d$, $N$, and, of course, the computational resources available and/or trade-off preferences between accuracy and speed. It also depends on the time step relative to the time scales of the flow. Second, to what extent does the improved convergence offered by the Toeplitz modes over the Hankel modes impact performance in applications of interest? Finally, to our knowledge, this study is the first to examine the convergence of space-time POD to space-only and spectral POD. By further analyzing the correlation matrix, it may be possible to make statements about when this convergence is to take place, e.g., that space-time POD performed on a window $T$ where $C(0,T) \approx 0$ will produce modes near to spectral POD modes.
\\

%%%%%%%%%%%%%%%%%%%%%%%%%%%%%%%%%%%%%%%%%%%%%%%%%%%%%%%%%%%%%%%%%%%%%%%%%%%%%%%
% -----------------------------------------------------------------------------
% -- Backmatter 
% -----------------------------------------------------------------------------

\section*{Acknowledgements}
We thank Ms. Peijing Liu for her contributions during the early stages of this project.
\bibliographystyle{siam}
\bibliography{references_abbrv_sansurl}

\end{document}